\documentclass[12pt,a4paper,leqno]{amsart}%{amsproc}
\usepackage{amsfonts,amssymb,amsmath,amstext,amscd,latexsym,euscript,amsthm}
\usepackage[T1]{fontenc}
\advance\textwidth by 3.2cm %%3.6
\advance\oddsidemargin by -1.6cm \advance\evensidemargin by -1.6cm
\input{colordvi}
\newtheoremstyle{plain}{2pt}{2pt}{\normalfont\sf}{\parindent}{\normalfont\bf}{.\,\,}{.0em}{}
\newtheoremstyle{definition}{2pt}{2pt}{\normalfont\sl}{\parindent}{\normalfont\bf}{.\,\,}{.0em}{}
\newtheoremstyle{proof}{2pt}{2pt}{\normalfont}{\parindent}{\normalfont\scshape}{.}{.0em}{}
\theoremstyle{definition}
\newtheorem{definition}{Definition}[section]

\newtheorem{remark}{Remark}[section]
\newtheorem*{definition*}{Definition}
\newtheorem*{example*}{Example}
\newtheorem*{remark*}{Remark}
\theoremstyle{proof}

\theoremstyle{plain}
\newtheorem{theorem}{Theorem}[section]
\newtheorem{lemma}{Lemma}[section]
\newtheorem{proposition}{Proposition}[section]
\newtheorem{corollary}{Corollary}[section]
\newtheorem{statement}{Statement}[section]

\newtheorem*{theorem*}{Theorem}
\newtheorem*{lemma*}{Lemma}
\newtheorem*{proposition*}{Proposition}
\newtheorem*{corollary*}{Corollary}
\newtheorem*{statement*}{Statement}
\newtheorem*{problem*}{Problem}
\newtheorem*{hypothesis*}{Hypothesis}
\DeclareMathOperator{\span1}{span}

\newcommand{\mmb}[1]{\mathbb{#1}}
\newcommand{\mmf}[1]{\mathbf{#1}}
\newcommand{\mmm}[1]{\mathrm{#1}}
\newcommand{\mcl}[1]{\mathcal{#1}}

\newcounter{num}

\newcounter{rlab}

\newcommand{\ov}[1]{\overline{#1}}

\def\Nat{{\mathbb N}}
\numberwithin{equation}{section}
\begin{document}
\title{Chernoff and Trotter type product formulas}
\author{Neklyudov A.}

\date{\today}

 \maketitle
\begin{abstract} We consider the abstract Cauchy problem
 $\dot{x}={\mmm{A}}x$, $x(0)=x_0\in\mcl{D}(\mmm{A})$ for linear
operators $\mmm{A}$ on a Banach space $\mmf{X}$. We prove uniqueness
of the (local) solution of this problem for a natural class of
operators $\mmm{A}$. Moreover, we establish that the solution
$x(\cdot)$ can be represented as a limit
$\lim\limits_{n\to\infty}\{\mmm{F}(t/n)^{n}\}$ in the weak operator
topology, where a function
$\mathrm{F}:[0,\infty)\mapsto\mathcal{L}(\mathbf{X})$ satisfies
$\mmm{F}'(0)y=\mmm{A}y$, $y\in \mcl{D}(\mmm{A})$. As a consequence,
we deduce necessary and sufficient conditions that a linear operator
$\mmm{C}$ is closable  and its closure is a generator of
$C_0$-semigroup. We also obtain some criteria for the sum of two
generators of $C_0$-semigroups to be  a generator of $C_0$-semigroup
such that the Trotter formula is valid.
\end{abstract}\maketitle
\tableofcontents

\section{Introduction.} Chernoff's Theorem  can be formulated as
follows (cf. \cite{Chern}, \cite{Engel}).

\begin{theorem*} Let $\mmm{F}$  be a map from $[0, \infty)$ to the space
of all continuous linear operators $\mcl{L}(\mmf{X})$ in a Banach
space $\mmf{X}$; $\mmm{Z}$ be a densely defined linear operator in
$\mmf{X}$ and the following conditions are satisfied:

\noindent i) $\mathrm{F}(0)=\mathrm{I}$ and there exists
$a\in\mathbb{R}$, $M\geq 1$ such that
$\|\mathrm{F}^{m}(\frac{s}{n})\|\leq M\exp( \frac{am}{n}s)$, for any
$n,\, m\in\mmb{N}$, $s\geq 0$.

\noindent ii) $\mcl{D}(\mmm{F}'(0))\supset\mcl{D}(\mmm{Z})$ and
$\mmm{F}'(0)f=\mmm{Z}f,\,\,f\in\mcl{D}(\mmm{Z})$.

\noindent iii) The closure of $\mmm{Z}$ exists and it is a generator
of $C_0$-semigroup.

\noindent Then  $\mmm{F}(t/n)^{n}$ converges to $\exp{(t
\overline{\mmm{C}})}$ as $n\to\infty$ in the strong operator
topology uniformly with respect to $t\in[0, T]$ for each $T>0$.
\end{theorem*}
 One of the applications of Chernoff's Theorem is the theory of the Feynman integrals (survey of this theory
can be found in the book \cite{SmoShav}). In the paper \cite{Hamil}
Chernoff's Theorem has been used to prove the representation of the
solution of the Schr\"{o}dinger equation as the Feynman integral.
Similar approach also has been used in  \cite{Weiz}, \cite{Witt}.

The main problem of application of Chernoff's Theorem is to check
condition (iii).
 In this work we study the class
$\EuScript{Z}$ of all densely defined operators $\mmm{Z}$ that
satisfy only conditions (i)--(ii) of Chernoff's Theorem. We can
notice that condition (iii) implies the existence of a function
$\mmm{F}$ satisfying conditions (i)--(ii) for a given $\mmm{Z}$.
Indeed, if condition (iii) is valid, then it is enough to put
$\mmm{F}(t)=\exp{(t\ov{\mmm{Z}})}$, $t\in[0,\,\infty)$. Thus, the
natural problem is to find
 conditions that a linear operator $\mmm{Z}\in\EuScript{Z}$ is closable and its closure is a generator of
$C_0$-semigroup. It turns out that conditions (i)--(ii) are
sufficient for closability of an operator $\mmm{Z}$ (see Theorem
\ref{cor:1}). Furthermore, if there exists (local) solution of the
abstract Cauchy problem $\dot{x}=\overline{\mmm{Z}}x$,
$x(0)=x_0\in\mcl{D}(\overline{\mmm{Z}})$, then this solution is
unique and can be represented as a limit
$\lim\limits_{n\to\infty}\{\mmm{F}(t/n)^{n}x_0\}$ in the weak
 topology $\sigma(\mmm{X},\,\mmm{X}^*)$ (see Theorem \ref{cor:2}). Applying these results, we
prove necessary and sufficient conditions that a linear operator
$\mmm{A}$ is closable and its closure is a generator of
$C_0$-semigroup (see Theorems \ref{cog:1}--\ref{t:54}, Corollaries
\ref{news:3}--\ref{new:1}). As a consequence, we deduce criteria for
the sum of two generators of $C_0$-semigroups to be  a generator of
$C_0$-semigroup such that the Trotter formula is valid (see
Corollaries \ref{t:3}--\ref{t:55}, \ref{corol}).

\section{Preliminaries.} For any normed space $(\mmf{E},\,\|\cdot\|_{\mmf{E}})$ (on the field
$\mmb{T}\in\{\mmb{R},\mmb{C}\}$) let $\mcl{L}(\mmf{E})$ be the space
of all bounded linear operators in $\mmf{E}$ with the topology of
pointwise convergence (strong operator topology), $\mmm{I}$ be the
identity operator in $\mmf{E}$. For any linear operator $\mmm{A}$ in
$\mmf{E}$ let $\mcl{D}(\mmm{A})$ be the domain of $\mmm{A}$. Linear
 subspaces of $\mmf{E}$ are considered as the normed spaces with the natural
norms (i. e. the norms inherited from the normed space $\mmf{E}$).
The dual space $\mmf{E}^{*}$ is the Banach space of all linear
continuous functionals on $\mmf{E}$ with the  norm
$\|\cdot\|_{\mmf{E}^*}$ defined by
$\|f\|_{\mmf{E}^*}=\sup\limits_{\|x\|_{\mmf{E}}\leq 1}{f(x)},\,
x\in\mmf{E},$ for each $f\in\mmf{E}^*$. Further, we will omit the
index of the space in the norm $\|\cdot\|_{\mmf{E}}$. For any
function $\mmm{S}:[0, \infty) \mapsto\mcl{L}(\mmf{E})$ we denote by
$\mmm{S}^*$ the function from $[0, \infty)$ to $\mcl{L}(\mmf{E})$
such that $\mmm{S}^*(s)=(\mmm{S}(s))^*$ for each $s\geq 0$.

\begin{definition} We call the (strong) derivative  at the point $0$ of a function
$\mmm{S}:[0, \infty) \mapsto\mcl{L}(\mmf{E})$ the linear map
$\mmm{S}'(0):\mcl{D}(\mmm{S}'(0))\mapsto\mmf{E}$ defined by
$\mmm{S}'(0)\psi=\lim\limits_{h\to
0}{h^{-1}(\mmm{S}(h)\psi-\mmm{S}(0)\psi)},
\,\psi\in\mcl{D}(\mmm{S}'(0)),$ where $\mcl{D}(\mmm{S}'(0))$ is the
space of all $\psi\in\mmf{E}$ such that the limit exists.
\end{definition}

\begin{definition}\label{r:1} The set $\EuScript{F}^{\mmf{E}}_{M,a}$, $M\geq 1$, $a\in\mmb{R}$, is the set of all
functions $\mathrm{F}:[0,\infty)\mapsto\mathcal{L}(\mathbf{E})$ that
satisfy the following conditions:

\noindent (i) $\mathrm{F}(0)=\mathrm{I}$ and
$\|\mathrm{F}^{m}(\frac{s}{n})\|\leq M\exp( \frac{am}{n}s)$ for all
$n,\, m\in\mmb{N}$, $s\geq 0$.

\noindent (ii) $\mcl{D}(\mmm{F}'(0))$ is dense in $\mathbf{E}$.
\end{definition}

\begin{definition} The set $\EuScript{F}_{\mmf{E}}$ is the set of all
functions $\mathrm{F}:[0,\infty)\mapsto\mathcal{L}(\mathbf{E})$ for
which there exists $M\geq 1$, $a\in\mmb{R}$ such that
$\mathrm{F}\in\EuScript{F}^{\mmf{E}}_{M,a}.$
\end{definition}

\begin{remark}
Condition (i) of Definition~\ref{r:1} can be replaced by
 equivalent condition (i*):

\noindent  (i*) $\mathrm{F}(0)=\mathrm{I}$ and
$\|\mathrm{F}^{n}(\frac{s}{n})\|\leq M\exp( as)$ for all
$n\,\in\mmb{N}$, $s\geq 0$.
\end{remark}

\begin{definition} The set $\EuScript{Z}$ is the set of all
densely defined linear operators in $\mathbf{E}$ for which there
exists $\mathrm{F}\in\EuScript{F}_{\mmf{E}}$ such that
$\mcl{D}(\mmm{Z})\subset\mcl{D}(\mmm{F}'(0))$ and
$\mmm{Z}f=\mmm{F}'(0)f$ for each $f\in\mcl{D}(\mmm{Z})$.
\end{definition}

\begin{definition}
 The function $\mmm{T}:[0,\,\infty)\to\mcl{L}(\mmf{E})$
is called $C_0$-semigroup if the following conditions are satisfied:

\noindent 1) $\mmm{T}(0)=\mmm{I}$,
$\mmm{T}(l+m)=\mmm{T}(l)\mmm{T}(m)$ for each $l,\,m\in[0,\,\infty)$.

\noindent 2) The function $\mmm{T}$ is continuous.

\noindent 3) There exists $M>1,\,\,a\in\mmb{R}$ such that
$\|\mmm{T}(s)\|\leq M \exp{(a s)}$ for each $s\geq 0$.
\end{definition}

\begin{definition}
The linear operator $\mmm{Z}$ is called the generator of
$C_0$-semigroup $\mmm{T}$ if $\mmm{Z}$ is the (strong) derivative at
the point $0$ of the function $\mmm{T}$.
\end{definition}
It is a well-known fact that there exists one-to-one correspondence
between $C_0$-semigroups and its generators.  The following results
can be found in \cite{Dav}:

\begin{statement}\label{st:31} Let $\mmf{E}$ be a Banach space. Assume that the
function $\mmm{T}:[0,\,\infty)\to\mcl{L}(\mmf{E})$ is
$C_0$-semigroup. Then $\mcl{D}({\mmm{T}^{*}}'(0))$ is *-dense in
$\mmf{E}^{*}$, where ${\mmm{T}^{*}}'(0)$ is the (strong) derivative
at the point $0$ of the function $\mmm{T}^{*}$.
\end{statement}
\begin{statement}\label{st:32} Let $\mmf{E}$ be a reflexive  Banach space. Assume that the
function $\mmm{T}:[0,\,\infty)\to\mcl{L}(\mmf{E})$ is
$C_0$-semigroup. Then the function $\mmm{T}^*$ is also
$C_0$-semigroup and its generator
${\mmm{T}^{*}}'(0)=({\mmm{T}}'(0))^{*}.$
\end{statement}
In the paper we assume that $\mmf{X}$
 is a Banach space and $\mcl{B}$ is a reflexive separable Banach space or a Hilbert
space.

\section{Closability of operators from class $\EuScript{Z}$.}

\begin{proposition}\label{t:1}
 Assume that  $\mathrm{F}\in\EuScript{F}_{\mmf{X}}$, $t>0$ and $\{n_k\}_{k=1}^{\infty}$ is an increasing sequence of natural
numbers.
 Then for any separable closed linear
spaces $\Phi\subset\mmf{X}$, $\Psi\subset\mmf{X}^*$ there exists a
subsequence $\{g_k\}_{k=1}^{\infty}$ of the sequence
$\{n_k\}_{k=1}^{\infty}$ such that there exists
\begin{eqnarray}
 \lim\limits_{k
\rightarrow\infty}{(\{\mathrm{F}(t/{g_k})\}^{[{g_k}s]}
g,\phi)}\label{eqn:l1}
\end{eqnarray}
for any $g\in \Phi$, $\phi\in \Psi$ and $s>0$. Furthermore, if such
subsequence $\{g_k\}_{k=1}^{\infty}$ is chosen, then the family of
the functions $\mathrm{T}_s: \Phi\times\Psi\mapsto\mathbb{T}$,
 $s\geq 0$, defined by
\begin{equation}
\mathrm{T}_s (g,\,\phi)=\lim\limits_{k
\rightarrow\infty}{(\{\mathrm{F}(t/{g_k})\}^{[{g_k}\frac{s}{t}]}
g,\phi)},\label{eqn:l2}
\end{equation}
  satisfies the following conditions:

\noindent a) $\mathrm{T}_s (g, \phi)$ is linear w.r.t. $g$ and
$\phi$ and the following inequality is satisfied $$\|\mathrm{T}_s
(g, \phi)\|\leq M \exp{(a s)}\|g\|\|\phi\|,\,\,s\geq 0.$$

\noindent b) $\mathrm{T}_s (g, \phi)$ is continuous w.r.t. $s$ for
any $g\in\Phi$, $\phi\in \Psi$ and the continuity is uniform w.r.t.
$\phi\in \Upsilon$, where $\Upsilon$ is a bounded subset of $\Psi$.

\noindent c) If $f\in\Phi\cap\mcl{D}(\mmm{F}'(0))$ and
$\mathrm{F}'(0)f\in\Phi$ then there exists $$(\mathrm{T}_{s} (f,
\phi))^{'}_{s}=\mathrm{T}_{s}(\mathrm{F}'(0)f, \phi)$$ for each
$s\geq 0$, $\phi\in \Psi$ uniformly w.r.t. $\phi\in \Upsilon$, where
$\Upsilon$ is a bounded subset of $\Psi$.

\noindent d) Assume $\phi\in\mcl{D}({\mmm{F}^*}'(0))$ and
${\mmm{F}^*}'(0)\phi\in\Psi$. Then $$\mathrm{T}_{m} (f,
\phi)-\mathrm{T}_{l} (f, \phi)=\int_{l}^{m}{\mathrm{T}_{s}(f,
{\mmm{F}^*}'(0)\phi)\,d s},\,\,m,\,\,l\geq 0,$$ for any $f\in\Phi$.

\noindent e) If $\phi\in\mcl{D}({\mmm{F}^*}'(0))$ and
${\mmm{F}^*}'(0)\phi\in\Psi$, then $$(\mathrm{T}_{s} (f,
\phi))^{'}_{s}=\mathrm{T}_{s}(f, {\mmm{F}^*}'(0)\phi),\,\,s\geq 0,$$
for any $f\in\Phi$.
\end{proposition}

Now we will formulate some auxiliary lemmas.

\begin{lemma}\label{t:5}
Assume that  $\mathrm{F}\in\EuScript{F}^{\mmf{X}}_{M,a}$, $M\geq 1$,
and $a=0$. Then there exists $\lim_{ i/k \to
0}{({\mathrm{F}}^{i}(t/k)-\mathrm{I})g}=0,$ $i, k \in \mathbb{N},$
for any $g \in \mmf{X}$.
\end{lemma}

\begin{remark} The limit in  Lemma~\ref{t:5} is considered as a limit
w.r.t. the filter base consisting of the sets of the form
$\{(i,k)|\,|\frac{i}{k}|<\epsilon,\,i,\, k \in \mathbb{N}\}$,
$\epsilon>0$.
\end{remark}

\begin{proof}[Proof of Lemma~\ref{t:5}.]
If $g \in \mcl{D}(\mmm{F}'(0))$, then the following chain of the
inequalities is satisfied:

\begin{eqnarray}
\|({\mathrm{F}}^{i}(t/k)-\mathrm{I})g\|&=& i
{\textstyle\frac{t}{k}}{\bigl\|\bigl({\textstyle\frac{\mathrm{F}(\frac{t}{k})^{i-1}
+\cdots+\mathrm{I}}{i}\bigr)
\frac{(\mathrm{F}(\frac{t}{k})g-g)}{\frac{t}{k}}}\bigr\|}\nonumber\\&\leq&
 i {\textstyle\frac{t}{k}}M\bigl\|{\textstyle\frac{(\mathrm{F}(\frac{t}{k})g-g)}{\frac{t}{k}}}\bigr\|\leq
i {\textstyle\frac{t}{k}}2 M \| \mmm{F}'(0)g \|,\nonumber
\end{eqnarray}

where $k>k_{0}(g)$. Therefore, we have $$\lim_{ i/k \to
0}{({\mathrm{F}}^{i}(t/k)-\mathrm{I})g}=0.$$

If $g \notin \mcl{D}(\mmm{F}'(0))$, then for any $\varepsilon >0$
there exists $g'$ such that $\|g-g'\|<\varepsilon /2 M$ and

\begin{eqnarray}
 \|\limsup\limits_{ i/k \to
0}{({\mathrm{F}}^{i}(t/k)-\mathrm{I})g}\|&\leq&
\|\limsup\limits_{ i/k \to 0}{({\mathrm{F}}^{i}(t/k)-\mathrm{I})(g-g')}\|\nonumber\\
 &+&\|\limsup\limits_{ i/k \to
0}{({\mathrm{F}}^{i}(t/k)-\mathrm{I})g'}\|\leq\varepsilon\nonumber
\end{eqnarray}

Since $\varepsilon$ is arbitrary the Lemma is proved.
\end{proof}

We can immediately deduce the following Lemma from Lemma~\ref{t:5}.

\begin{lemma}\label{t:7}
Assume that $\mathrm{F}\in\EuScript{F}^{\mmf{X}}_{M,a}$, $M\geq 1$,
and $a=0$. Then for any $g \in \mmf{X}$ there exists $\lim\limits_{
|i-l/k| \to 0}{({\mathrm{F}}^{i}(t/k)-{\mathrm{F}}^{l}(t/k))g}=0 $,
$i, l, k \in \mathbb{N}$.
\end{lemma}
\begin{remark} The limit in  Lemma~\ref{t:7} is considered as a limit
w.r.t. the filter base consisting of the sets of the form
$\{(i,l,k)|\,0<|\frac{i-l}{k}|<\epsilon,\,i,\, k,\, l \in
\mathbb{N}\}$, $\epsilon>0$. In the same sense we will consider the
limit in  Lemma~\ref{t:8}.
\end{remark}

\begin{lemma}\label{t:8}
Assume that $\mathrm{F}\in\EuScript{F}^{\mmf{X}}_{M,a}$, $M\geq 1$,
and $a=0$. Then for any $g \in \mcl{D}(\mmm{F}'(0))$ there exists
$$\lim_{ |i-l/k| \to
0}{\frac{({\mathrm{F}}^{i}(t/k)-{\mathrm{F}}^{l}(t/k))g}{(i-l)t/k}-
({\mathrm{F}}^{\min(i,l)}(t/k))\mmm{F}'(0)g}=0,$$
 $i, l, k \in \mathbb{N}$.
\end{lemma}

\begin{proof}[Proof of Lemma~\ref{t:8}.] We have the following chain of the inequalities

\begin{eqnarray}
&\,&\limsup\limits_{ |i-l/k| \to
0}{\|{\textstyle{\frac{({\mathrm{F}}^{i}(\frac{t}{k})-{\mathrm{F}}^{l}(\frac{t}{k})g}{(i-l)t/k}}-
({\mathrm{F}}^{\min(i,l)}(\frac{t}{k}))\mmm{F}'(0)g}\|}\nonumber\\
&\leq& M\limsup\limits_{ |i-l/k| \to
0}{\|{\textstyle\frac{({\mathrm{F}}^{|i-l|}(\frac{t}{k})-\mathrm{I})g}{|i-l|\frac{t}{k}}}-\mmm{F}'(0)g}\|
 \leq M\limsup\limits_{ i/k \to
0}{\|{\textstyle\frac{({\mathrm{F}}^{i}(\frac{t}{k})-\mathrm{I})g}{i
\frac{t}{k}}}-\mmm{F}'(0)g}\|\nonumber \\ &\leq& M\limsup\limits_{
i/k \to
0}{{\bigl\|\bigl({\textstyle\frac{\mathrm{F}(\frac{t}{k})^{i-1}
+\cdots+\mathrm{I}}{i}}\bigr) ({\textstyle
\frac{(\mathrm{F}(\frac{t}{k})g-g)}{\frac{t}{k}}}-\mmm{F}'(0)g)}\bigr\|}\nonumber
\\
 &+& M\limsup\limits_{ i/k \to
0}{{\bigl\|\bigl({\textstyle\frac{\mathrm{F}(\frac{t}{k})^{i-1}
+\cdots+\mathrm{I}}{i}}\bigr)\mmm{F}'(0)g
-\mathrm{F}'(0)g}\bigr\|}\nonumber \\
&\leq&
 \limsup\limits_{ i/k \to 0}{M^2{\bigl\|
({\textstyle\frac{(\mathrm{F}(\frac{t}{k})g-g)}{\frac{t}{k}}}-\mmm{F}'(0)g)}\bigr\|}\nonumber
\\
 &+& M\limsup\limits_{ i/k \to
0}{{\bigl\|\bigl({\textstyle\frac{(\mathrm{F}(\frac{t}{k})^{i-1}-\mathrm{I})
+\cdots+(\mathrm{I}-\mathrm{I})}{i}}\bigr)\mmm{F}'(0)g}\bigr\|}=0,\nonumber
\end{eqnarray}
where the last limit is equal to $0$ by Lemma $1$.
\end{proof}

\begin{proof}[Proof of Proposition~\ref{t:1}.]
We consider only the case $a = 0$. The case  $a \neq 0$ can be
easily reduced to the case $a=0$ by examining $\mathrm{F}(s)\exp(-a
s)$ instead of $\mathrm{F}(s)$.

We can choose a subsequence $\{g_k\}_{k=1}^{\infty}$ of the sequence
$\{n_k\}_{k=1}^{\infty}$ such that for any nonnegative rational $r$
there exists
$$\lim_{k\to\infty}{(\{\mathrm{F}(t/{g_k})\}^{[{g_k}r]} g,\phi)},
\,\,g\in\mmf{C},\,\,\phi\in \mmf{D},$$ where $\mmf{C}$ and $\mmf{D}$
are dense countable subsets of $\Phi$ and $\Psi$ respectively.
Indeed, it is possible to find such subsequence because
$$\|\{\mathrm{F}(t/{n_k})\}^{[{n_k}r]}\|\leq M$$ and we have a
countable number of conditions on the sequence
$\{g_k\}_{k=1}^{\infty}$. Hence, it follows from the density of
$\mmf{C}$ in $\Phi$ (and the density of $\mmf{D}$ in $\Psi$) that
there exists
$$\lim_{k\to\infty}{(\{\mathrm{F}(t/{g_k})\}^{[{g_k}r]} g,\phi)}$$
for any $g\in\Phi$, $\phi\in\Psi$. Let us show that the last limit
exists for any real number $r\geq 0$. Fix $r>0$, $z\in\Psi$ and
$\epsilon
> 0$. Put $\epsilon^{'}=\epsilon/3\|z\|$. We can choose
$\epsilon^{''}(g,\epsilon^{'})$ such that
$$\|({\mathrm{F}}^{i}(t/k)-{\mathrm{F}}^{l}(t/k))g\|<\epsilon^{'},$$
when $|(i-l)/k|<\epsilon^{''}$,$i,l,k\in\mmb{N}$. Indeed, the
existence of $\epsilon^{''}$ follows from Lemma~\ref{t:7}. Choose
also a positive $s\in\mmb{Q}$ such that $|r-s|<\epsilon^{''}/2$.
There exists $k_0\in\mmb{N}$ such that $2/g_{k_0}<\epsilon^{''}/2$.
Now it follows from the existence of the limit
$$\lim\limits_{k\to\infty}\{(\mathrm{F}(t/{g_k})\}^{[{g_k}s]}g,z)$$ that we can find $n_0>k_0$ such that for any $l,m>
n_0$, $l,m\in\mmb{N},$ we have the following inequality
$$|((\{\mathrm{F}(t/{g_l})\}^{[{g_l}s]}-
\{\mathrm{F}(t/{g_m})\}^{[{g_m}s]})g, z)|<\epsilon/3.$$ Therefore,
we get the following chain of the inequalities

\begin{eqnarray}
 &\,&|((\{\mathrm{F}(t/{g_l})\}^{[{g_l}r]}-
\{\mathrm{F}(t/{g_m})\}^{[{g_m}r]})g, z)|\leq
|((\{\mathrm{F}(t/{g_l})\}^{[{g_l}r]}\nonumber\\&-&
\{\mathrm{F}(t/{g_l})\}^{[{g_l}s]})g, z)|
 +|((\{\mathrm{F}(t/{g_l})\}^{[{g_l}s]}-
\{\mathrm{F}(t/{g_m})\}^{[{g_m}s]})g, z)|\nonumber\\&+&
|((\{\mathrm{F}(t/{g_m})\}^{[{g_m}s]}-
\{\mathrm{F}(t/{g_m})\}^{[{g_m}r]})g, z)|
 \leq \epsilon/3 +\epsilon/3+\epsilon/3\leq\epsilon,\nonumber
\end{eqnarray}
where we have used the fact that $$|[g_n r]-[g_n s]|/g_n \leq
|g_nr-g_n s|+2/g_n \leq\epsilon^{''},\,\,n>k_0.$$ As a consequence,
we get that there exists
$$\lim_{k\to\infty}{(\{\mathrm{F}(t/{g_k})\}^{[{g_k}s]}
g,\phi)},\,\, s\geq 0,\,\, g\in\Phi,\,\,\phi\in\Psi.$$ Thus, the
first part of Proposition~\ref{t:1} is  proved.

Now let us consider parts (a) -- (e) of Proposition~\ref{t:1}.

\begin{trivlist}
\item[(a)] Part $(a)$ immediately follows from the definition of
$\mmm{T}_s$. \item[(b)] Fix $g\in\Phi$. It follows from
Lemma~\ref{t:7} that for any $\epsilon>0$ there exist $h_0>0$,
$k_0\in\mmb{N}$ such that
$$|(\{\mathrm{F}(t/{g_k})\}^{[{g_k}\frac{s+h}{t}]} g,\phi)-
(\{\mathrm{F}(t/{g_k})\}^{[{g_k}\frac{s}{t}]}
g,\phi)|<\epsilon\|\phi\|,$$ for any $|h|<h_0,\,h\geq -s$,
$k>k_0,\,k\in\mmb{N}$. Letting $k\to\infty$, we get the inequality
$$|\mmm{T}_{s+h}(g,\phi)-\mmm{T}_{s}(g,\phi)|\leq\epsilon\|\phi\|.$$
From the arbitrariness of $\epsilon$ we easily infer part $(b)$ of
Proposition~\ref{t:1}.

\item[(c)] Assume that $f\in\Phi\cap\mcl{D}(\mmm{F}'(0))$,
$\mathrm{F}'(0)f\in\Phi$. It follows from Lemma~\ref{t:8} that for
any $\epsilon>0$ there exist $h_0>0$, $k_0\in\mmb{N}$ such that
\begin{eqnarray}
\nonumber & & |({\mathrm{F}}^{[{g_k}\frac{s}{t}]}(t/g_k)f,\phi)-
 ({\mathrm{F}}^{
[{g_k}\frac{s+h}{t}]}(t/g_k)f,\phi)-t/g_k([{g_k}{\textstyle\frac{s}{t}}]-
[{g_k}{\textstyle\frac{s+h}{t}}])\\&\times&
({\mathrm{F}}^{\min([{g_k}\frac{s}{t}],[{g_k}\frac{s+h}{t}])}(t/g_k)\mathrm{F}'(0)f,\phi)|\nonumber
<\epsilon\|\phi\||[{g_k}{\textstyle\frac{s}{t}}]-[{g_k}{\textstyle\frac{s+h}{t}}]|{\textstyle\frac{t}{g_k}}
\end{eqnarray}
for any $|h|<h_0,\,h\geq -s$, $k>k_0,\,k\in\mmb{N}$. Letting
$k\to\infty$, we get the inequality
\begin{equation}
\|\mathrm{T}_{s}(f,\phi)-\mathrm{T}_{s+h}(f,\phi)-h\mathrm{T}_{\min{(s,
s+h)}}(\mathrm{F}'(0)f,\phi)\|\leq \epsilon |h|\|\phi\|. \nonumber
\end{equation}
From the arbitrariness of $\epsilon$ we deduce part $(c)$ of
Proposition~\ref{t:1}.

\item[(d)]
 Assume that $\phi\in\mcl{D}({\mmm{F}^*}'(0))$,
${\mmm{F}^*}'(0)\phi\in\Psi$. Fix $\epsilon>0$ and $m>l$. Choose
$r_0$, $j_0$ such that
$$\|{\mathrm{F}(t/{g_k})}^{[{g_k}\frac{g}{t}]})f-{\mathrm{F}(t/{g_k})}^{[{g_k}\frac{h}{t}]})f\|
<\epsilon, \,\,k\geq j_0, \,\,|g-h|<r_0,\,\,0\leq g,h<m.$$ The
existence of such $r_0$, $j_0$ follows from Lemma~\ref{t:7}. Hence,
we have
\begin{eqnarray}|\mmm{T}_h(f,{\mathrm{F}^{*}}'(0)\phi)-\mmm{T}_g(f,{\mathrm{F}^{*}}'(0)\phi)|\leq\epsilon,\label{q:1}
\end{eqnarray} where $|g-h|<r_0,\,0\leq g,h<m$.
Therefore, we get
\begin{eqnarray}
& &|\sum_{v=0}^{v_0-1}{\mathrm{T}_{l+v(m-l)/v_0}(f,
{\mmm{F}^*}'(0)\phi)((m-l)/v_0)}\nonumber\\&-&\int_{l}^{m}{\mathrm{T}_{s}(f,
{\mmm{F}^*}'(0)\phi)\,d s}|\leq \epsilon |m-l|,\label{eqn:d0}
\end{eqnarray}
$v_0=[(m-l)/r_0]+1$.
 Choose $j_1>j_0$ such that
\begin{eqnarray}
|\mathrm{T}_{l+v(m-l)/v_0}(f, {\mmm{F}^*}'(0)\phi)-
({\mathrm{F}(t/{g_k})}^{[{g_k}(l+v(m-l)/v_0)/t]}f,{\mathrm{F}^{*}}'(0)\phi)|<\epsilon\nonumber
\end{eqnarray}
for any $v\in\{0,1,\dots ,v_0\}$, $k>j_1$. Thus,
\begin{eqnarray}
&\,&|\sum_{v=0}^{v_0-1}{\mathrm{T}_{l+v(m-l)/v_0}(f,
{\mmm{F}^*}'(0)\phi)((m-l)/v_0)}\label{eqn:d1}\\
&-&((m-l)/v_0)\sum_{v=0}^{v_0-1}{({\mathrm{F}(t/{g_k})}
^{[{g_k}(l+v(m-l)/v_0)/t]}f,{\mathrm{F}^{*}}'(0)\phi)}|<\epsilon
|m-l|.\label{eqn:d2}
\end{eqnarray}
Now we notice that
\begin{eqnarray}
&\,&\mathrm{T}_{m}(f,\phi)-\mathrm{T}_{l}(f,\phi)=\lim\limits_{k
\rightarrow\infty}{({\mathrm{F}({\textstyle\frac{t}{g_k}})}^{[{g_k}\frac{m}{t}]}
f-{\mathrm{F}({\textstyle\frac{t}{g_k}})}^{[{g_k}\frac{l}{t}]} f,\phi)}\label{q:2}\\
&=&\lim\limits_{k
\rightarrow\infty}{((\mathrm{F}({\textstyle\frac{t}{g_k}})
-\mathrm{I})
({\mathrm{F}({\textstyle\frac{t}{g_k}})}^{[{g_k}\frac{m}{t}]-1}+\dots+{\mathrm{F}({\textstyle\frac{t}{g_k}})}^{[{g_k}\frac{l}{t}]})f,\phi)}\label{q:3}\\
&=&\lim\limits_{k \rightarrow\infty}{
({\textstyle\frac{t}{g_k}}({\mathrm{F}({\textstyle\frac{t}{g_k}})}^{[{g_k}\frac{m}{t}]-1}+\dots+{\mathrm{F}({\textstyle\frac{t}{g_k}})}^{[{g_k}\frac{l}{t}]})f,
{\textstyle\frac{t}{g_k}}(\mathrm{F}^{*}({\textstyle\frac{t}{g_k}}) -\mathrm{I})\phi)}\label{q:4}\\
&=&\lim\limits_{k \rightarrow\infty}{
({\textstyle\frac{t}{g_k}}({\mathrm{F}({\textstyle\frac{t}{g_k}})}^{[{g_k}\frac{m}{t}]-1}+\dots+{\mathrm{F}({\textstyle\frac{t}{g_k}})}^{[{g_k}\frac{l}{t}]})f,
{\mathrm{F}^{*}}'(0)\phi)}\label{q:7}\\
&=&\lim\limits_{k
\rightarrow\infty}{{\textstyle\frac{t}{g_k}}\sum_{s=0}^{d_k}{({\mathrm{F}({\textstyle\frac{t}{g_k}})}^{[{g_k}(l+s({\textstyle\frac{t}{g_k}}))/t]}f,{\mathrm{F}^{*}}'(0)\phi)}}\label{q:8}\\
&=&\lim\limits_{k
\rightarrow\infty}{\sum_{v=0}^{v_0-1}{\sum_{s\in\mmf{B}_{k,v,v_0}}{{\textstyle\frac{t}{g_k}}({\mathrm{F}({\textstyle\frac{t}{g_k}})}^{[{g_k}(l+s({\textstyle\frac{t}{g_k}}))/t]}f,{\mathrm{F}^{*}}'(0)\phi)}}},\label{q:9}
\end{eqnarray}
where $d_k=[(m-l)g_k/t]$ and
$$\mmf{B}_{k,v,v_0}=\{s\in\mmb{N}|\,v(m-l)/v_0\leq s t/g_k <
(v+1)(m-l)/v_0\}.$$

\noindent So, by inequality~\eqref{q:1} and equalities \eqref{q:2}
--- \eqref{q:9}, we get
\begin{eqnarray}
|\mathrm{T}_{l}(f,\phi)&+&
\sum_{v=0}^{v_0-1}{({\mathrm{F}({\textstyle\frac{t}{g_k}})}
^{[{g_k}(l+v{\textstyle\frac{m-l}{v_0}})/t]}f,{\mathrm{F}^{*}}'(0)\phi)({\textstyle\frac{m-l}{v_0}})}-\mathrm{T}_{m}(f,\phi)| \label{eqn:d3}\\
&=&|\lim\limits_{k
\rightarrow\infty}{\sum_{v=0}^{v_0-1}{\sum_{s\in\mmf{B}_
{k,v,v_0}}{{\textstyle\frac{t}{g_k}}({\mathrm{F}({\textstyle\frac{t}{g_k}})}^{[{g_k}(l+s({\textstyle\frac{t}{g_k}}))/t]}f,{\mathrm{F}^{*}}'(0)\phi)}}}\\
&-&({\textstyle\frac{m-l}{v_0}})\sum_{v=0}^{v_0-1}{({\mathrm{F}({\textstyle\frac{t}{g_k}})}
^{[{g_k}(l+v{\textstyle\frac{m-l}{v_0}}/t]}f,{\mathrm{F}^{*}}'(0)\phi)}|\\
&\leq&\sum_{v=0}^{v_0-1}|\lim\limits_{k
\rightarrow\infty}{\sum_{s\in\mmf{B}_{k,v,v_0}}{{\textstyle\frac{t}{g_k}}({\mathrm{F}({\textstyle\frac{t}{g_k}})}^{[g_k(l+s({\textstyle\frac{t}{g_k}}))/t]}f,{\mathrm{F}^{*}}'(0)\phi)}}\\
&-&
{\textstyle\frac{m-l}{v_0}}(\mathrm{F}({\textstyle\frac{t}{g_k}})^{[g_k(l+v{\textstyle\frac{m-l}{v_0}})/t]}f,{\mathrm{F}^*}'(0)\phi)|<\epsilon
v_0 {\textstyle\frac{|m-l|}{v_0}}.\label{eqn:d4}
\end{eqnarray}
Combining inequalities \eqref{eqn:d0}, \eqref{eqn:d1} ---
\eqref{eqn:d2} and \eqref{eqn:d3}
--- \eqref{eqn:d4}, we get
\begin{eqnarray}
|\mathrm{T}_{m}(f,\phi)-\mathrm{T}_{l}(f,\phi)-\int_{l}^{m}{\mathrm{T}_{s}(f,
{\mmm{F}^*}'(0)\phi)\,d s}|<3\epsilon |m-l|\nonumber
\end{eqnarray}
From the arbitrariness of $\epsilon$ we deduce part $(d)$ of
Proposition~\ref{t:1}. \item[(e)] Part $(e)$ is a direct consequence
of parts $(b)$ and $(d)$.
\end{trivlist}
\end{proof}

\begin{proposition}\label{n:1}Assume that
$\mathrm{F}\in\EuScript{F}_{\mmf{X}}$. Then the operator
$\mmm{F}'(0)$ is closable.
\end{proposition}

\begin{proof}[Proof of Proposition~\ref{n:1}.] It is enough to show that the existence of the limits
$$\lim_{n\to\infty}{f_n}=0,\,\, \lim_{n\to\infty}{\mmm{F}'(0)f_n}=h,\,\, h\in\mmf{X},$$
for any sequence $\{ f_n \}_{n=1}^{\infty}$,
$f_n\in\mcl{D}(\mmm{F}'(0)),$ imply that $h=0.$ Indeed, it
immediately follows from Lemma 1.8, p. 7 in \cite{Dav}. We will
argue by contradiction. Assume that $h\neq 0$. Then there exists
$\phi\in\mmf{X}^{*}$ such that $(h,\,\phi)\neq 0$. Let $\mmf{B}$ be
the closure of the linear span of the elements
$f_i$,$\mmm{F}'(0)f_i$, $i\in\mmb{N}$.
 Then, according to Proposition~\ref{t:1}, there exists an increasing sequence $\{g_k\}_{k=1}^{\infty},\,\,g_k\in\mmb{N},$
 such that we can define the family of the functions
$\mathrm{T}_s: \Phi\times\Psi\mapsto\mathbb{T}$,
 $s\geq 0$, by formula \eqref{eqn:l2} with $\Phi=\mmf{B}$ and $\Psi=\{\phi\}$.
Hence, we have $$(\mathrm{T}_{s}
(f_i,\phi))^{'}_{s}=\mathrm{T}_{s}(\mathrm{F}'(0)f_i,\phi)$$ for any
$i\in\mmb{N}$. Therefore, we get
\begin{equation}
\mathrm{T}_{s}(f_i,\phi)-(f_i,\phi)=\int_0^{s}{\mathrm{T}_{t}(\mathrm{F}'(0)f_i,\phi)\,
d t},i\in\mmb{N}.\label{eqn:eq-2}
\end{equation}
Letting  $i\to \infty$ in expression \eqref{eqn:eq-2}, we get the
equality
$$
\mathrm{T}_{s}(0,\phi)-(0,\phi)=\int_0^{s}{\mathrm{T}_{r}(h,\phi)\,
d r}=0,s\geq 0.
$$
Indeed, it immediately follows from \eqref{eqn:eq-2} and the uniform
continuity of the family $\{\mathrm{T}_{l}(\cdot,\phi)|\, l\in [0,\,
s]\}$. Hence, we get
\begin{equation}
(\int_0^{s}{\mathrm{T}_{t}(h,\phi)\,
dt})'_s=\mathrm{T}_{s}(h,\phi)=0,\,\,s\geq 0.\label{eqn:eq-3}
\end{equation}
Put $s=0$ in \eqref{eqn:eq-3}. Consequently, we have
$$
\mathrm{T}_{s}(h,\phi)=(h,\phi)=0.
$$
Thus, we get contradiction with the assumption $h\neq 0$.
\end{proof}

As a consequence of Proposition~\ref{n:1}, we get
\begin{theorem}\label{cor:1}Let
$\mathrm{F}\in\EuScript{F}_{\mmf{X}}$. Assume that a linear operator
$\mmm{Z}$ has the domain
$\mcl{D}(\mmm{Z})\subset\mcl{D}(\mmm{F}'(0))$ and
$$\mmm{Z}f=\mmm{F}'(0)f,$$ for any $f\in\mcl{D}(\mmm{Z}).$ Then
$\mmm{Z}$ is closable.
\end{theorem}

\section{Chernoff type formula for solution of the abstract Cauchy problem.}

\begin{lemma}\label{lem:lem-4}
Suppose that the assumptions of Proposition~\ref{t:1} are satisfied.
Assume also that $f_i\in\mcl{D}(\mmm{F}'(0))\cap\Phi$,
$\mmm{F}'(0)f_i\in\Phi$ for any $i\in\mmb{N}$. Furthermore, suppose
that
$$
\lim_{i\to\infty}{f_i}=h,\,\, \lim_{i\to\infty}{\mmm{Z}f_i}=g.
$$
Then
$$(\mathrm{T}_{s}(h,\phi))^{'}_{s}=\mathrm{T}_{s}(g,\phi),\,\,s\geq 0,\,\,\phi\in\Psi.$$
\end{lemma}
\begin{proof}[Proof of Lemma~\ref{lem:lem-4}.]
Notice that in the equalities
$$
\lim_{l\to
0}{{l}^{-1}(\mmm{T}_{s+l}(f_i,\phi)-\mmm{T}_{s}(f_i,\phi))}=\lim_{l\to
0}{{l}^{-1}{\int_{s}^{s+l}{\mmm{T}_{m}(\mmm{Z}f_i,\phi)\, d
m}}}=\mmm{T}_{s}(\mmm{Z}f_i,\phi)
$$
we have the uniform convergence w.r.t. $i\in\mmb{N}$ because there
exists $$\lim_{i\to\infty}{\mmm{Z}f_i}=h$$ and the family of the
functions $\{\mathrm{T}_{l}(\cdot,\phi)|\, l\in [0,\, s]\}$ is
uniformly continuous. Therefore,

\begin{eqnarray}
\mmm{T}_s (g,\phi)&=&\lim_{i\to\infty}{\lim_{l\to
0}{{l}^{-1}(\mmm{T}_{s+l}(f_i,\phi)-\mmm{T}_{s}(f_i,\phi))}}\nonumber\\&=&\lim_{l\to
0}{\lim_{i\to\infty}{{l}^{-1}(\mmm{T}_{s+l}(f_i,\phi)-\mmm{T}_{s}(f_i,\phi))}}\nonumber\\
&=&\lim_{l\to
0}{{l}^{-1}(\mmm{T}_{s+l}(h,\phi)-\mmm{T}_{s}(h,\phi))}=(\mathrm{T}_{s}
(h,\phi))^{'}_{s}.\nonumber
\end{eqnarray}
\end{proof}
\begin{lemma}\label{n:2}Suppose that the assumptions of Proposition~\ref{t:1} are satisfied. If
$f\in\mcl{D}(\overline{\mmm{F}'(0)})\cap\Phi$ and
$\overline{\mmm{F}'(0)}f\in\Phi$, then
$$
(\mathrm{T}_{s}
(f,\phi))^{'}_{s}=\mathrm{T}_{s}(\overline{\mmm{F}'(0)}f,\phi),\,\,s\geq
0,\,\, \phi\in\Psi.
$$
\end{lemma}
\begin{proof}[Proof of Lemma~\ref{n:2}.]
Since $f\in\mcl{D}(\overline{\mmm{F}'(0)})$, we see that there
exists a sequence $\{f_i\}_{i=1}^{\infty}$,
$f_i\in\mcl{D}(\mmm{F}'(0)),$ such that
$\lim\limits_{i\to\infty}{f_i}=f$ and
$$\lim_{i\to\infty}{\mmm{F}'(0)f_i}=\overline{\mmm{F}'(0)}f.$$ Let
$\Phi'$ be the minimal closed space such that $\Phi\subset\Phi'$ and
$f_i,\, \mmm{F}'(0)f_i\in\Phi'$ for all $i\in\mmb{N}$. Then,
according to Proposition~\ref{t:1}, we can choose a subsequence
$\{m_k\}_{k=1}^{\infty}$ of the sequence $\{g_k\}_{k=1}^{\infty}$
such that we can define the family of the functions
$\mathrm{T}^1_{s}:(\Phi',\Psi)\mapsto\mmb{T},$ $s\geq 0,$ by the
equality
$$\mathrm{T}^1_{s}(g,\phi)=\lim\limits_{k
\rightarrow\infty}{(\{\mathrm{F}(t/{m_k})\}^{[{m_k}\frac{s}{t}]}
g,\phi)}, \,\,g\in\Phi', \,\,\phi\in\Psi.$$ Thus, it follows from
Lemma~\ref{lem:lem-4} (applied with parameters $h=f$,
$g=\overline{\mmm{F}'(0)}f$) that $$(\mathrm{T}^1_{s}
(f,\phi))^{'}_{s}=\mathrm{T}^1_{s}(\overline{\mmm{F}'(0)}f,\phi),\,\,s\geq
0.$$ Since the restriction of $\mathrm{T}^1_{s}$  to the set
$\Phi\times\Psi$ is equal to $\mathrm{T}_{s}$ we have
$$(\mathrm{T}_{s}
(f,\phi))^{'}_{s}=\mathrm{T}_{s}(\overline{\mmm{F}'(0)}f,\phi),\,\,s\geq
0,\,\,\phi\in\Psi,
$$
and the Lemma is proved.
\end{proof}
\begin{proposition}\label{t:6}Suppose that  $\mathrm{F}\in\EuScript{F}_{\mmf{X}}$ and $t,\,l>0$.
Assume also that there exists a local solution
$f:[0,\,l)\mapsto\mcl{D}(\overline{\mmm{F}'(0)})$ of the system
\begin{eqnarray}
f'(s) &=& \overline{\mmm{F}'(0)}f(s),\,\,s\in[0,\,l),\label{eqn:eq-4}\\
f(0) &=& g\in\mcl{D}(\overline{\mmm{F}'(0)}). \label{eqn:eq-41}
\end{eqnarray}
Then this solution is unique and
\begin{eqnarray}
f(s)=w\mbox{-}\lim_{n\to\infty}{\{\mathrm{F}(t/{n})\}^{[n\frac{s}{t}]}g},\,\,
s\in [0,\,l).\label{eqn:l3}
\end{eqnarray}
\end{proposition}
\begin{proof}[Proof of Proposition~\ref{t:6}.]Let $f(s),\,\,s\in [0,\,l),$ be a local solution of system
\eqref{eqn:eq-4}-\eqref{eqn:eq-41}. Then for any $r\in\mmb{Q}\cap
[0,\, l)$ there exists a sequence
$\{x_n^r\}_{n=1}^{\infty},\,\,x_n^r\in\mcl{D}(\mmm{F}'(0)),$ such
that
$$\lim_{n\to\infty}{x_n^r}=f(r)$$ and
$$\lim_{n\to\infty}{\overline{\mmm{F}'(0)}x_n^r}=\overline{\mmm{F}'(0)}f(r).$$
Let $\mmf{B}$ be the minimal closed linear subspace of $\mmf{X}$
such that $x_n^r\in\mmf{B}$ for all $r\in\mmb{Q}\cap [0,\, l),\,
n\in\mmb{N}$. Hence,  $f(s)\in\mmf{B}$ for all $s\in [0,\, l)$ and,
consequently, $\overline{\mmm{F}'(0)}f(s)\in\mmf{B}$ for all $s\in
[0,\, l)$. For any $\phi\in\mmf{X}^{*}$ let us choose a subsequence
$\{g_k\}_{k=1}^{\infty}$ of an arbitrary increasing sequence
$\{n_k\}_{k=1}^{\infty},\,\,n_k\in\mmb{N},$ such that we can define
the family of the functions $\mathrm{T}_s:
\Phi\times\Psi\mapsto\mathbb{T}$, $s\geq 0,$ by equality
\eqref{eqn:l2} with $\Phi=\mmf{B}$, $\Psi=\{\phi\}$. From
Lemma~\ref{n:2} it follows that
$$(\mmm{T}_s
(f(v),\phi))'_s=\mmm{T}_s(\overline{\mmm{F}'(0)}f(v),\phi),\,\,s,\,v\in[0,\,l).$$
  So, by Proposition~\ref{t:1},
\begin{eqnarray}
(\mmm{T}_s (f(v-s),\phi))'_s &=& \lim_{a\to
0}{{a}^{-1}(\mmm{T}_{s+a}(f(v-a-s),\phi)-\mmm{T}_{s}(f(v-s),\phi))}\nonumber\\&=&
\lim_{a\to
0}{{a}^{-1}(\mmm{T}_{s+a}((f(v-a-s)-f(v-s)),\phi))}\nonumber\\&+&
\lim_{a\to
0}{{a}^{-1}(\mmm{T}_{s+a}(f(v-s),\phi)-\mmm{T}_{s}(f(v-s),\phi))}\nonumber\\&=&
\mmm{T}_{s}(\overline{\mmm{F}'(0)}f(v-s),\phi))-\mmm{T}_{s}(\overline{\mmm{F}'(0)}f(v-s),\phi)=0,\nonumber
\end{eqnarray}
for $v\in (0,\,l),\,\,s\in(0,\,v).$ As a consequence, the function
$\mmm{T}_s (f(v-s),\phi)$ is a constant w.r.t. $s$. Thus,
$$\mmm{T}_v (g,\phi)=( f(v),\phi),\,\,v\in[0,\,l).$$ Now the
existence of limit~\eqref{eqn:l3} follows from the arbitrariness of
the sequence $\{n_k\}_{k=1}^{\infty}$ for any fixed
$\phi\in\mmf{X}^*$. Therefore, $f$ is unique.
\end{proof}

As a consequence of Proposition~\ref{t:6}, we get
\begin{theorem}\label{cor:2}Assume that $\mathrm{F}\in\EuScript{F}_{\mmf{X}}$ and $t,\,\,l>0$. Assume also that a linear operator $\mmm{Z}$
has the domain $\mcl{D}(\mmm{Z})\subset\mcl{D}(\mmm{F}'(0))$ and
$$\mmm{Z}f=\mmm{F}'(0)f,\,\,f\in\mcl{D}(\mmm{Z}).$$ If there exists
a local solution $f:[0,\,l)\mapsto\mcl{D}(\overline{\mmm{Z}})$ of
the system
\begin{eqnarray}
f'(s)&=&\overline{\mmm{Z}}f(s),\,\,s\in [0,\,l),\nonumber\\
f(0)&=&g\in\mcl{D}(\overline{\mmm{Z}}),\nonumber
\end{eqnarray}
then the solution is unique and
\begin{eqnarray}
f(s)=w\mbox{-}\lim_{n\to\infty}{\{\mathrm{F}(t/{n})\}^{[n\frac{s}{t}]}g},\,\,s\in
[0,\,l ).\nonumber
\end{eqnarray}
\end{theorem}
\begin{remark}
 The existence of
the closure of  $\mmm{Z}$ in Theorem~\ref{cor:2} follows from
Theorem~\ref{cor:1}.
\end{remark}
\begin{corollary}\label{cor_t:7}Assume that the conditions of Theorem~\ref{cor:2} are satisfied.
Then
$$f'(s)=w\mbox{-}\lim_{n\to\infty}{\{\mathrm{F}(t/{n})\}^{[n\frac{s}{t}]}
\overline{\mmm{Z}}f(0)},\,\,s\in[0,\, l).$$
\end{corollary}
\begin{proof}[Proof of Corollary~\ref{cor_t:7}.]
Define the set $\mmf{B}$  as in the proof of Proposition~\ref{t:6}.
For any $\phi\in\mmf{X}^{*}$ choose a subsequence
$\{g_k\}_{k=1}^{\infty}$ of an arbitrary increasing sequence
$\{n_k\}_{k=1}^{\infty},\,\,n_k\in\mmb{N},$ such that we can define
the family of the functions $\mathrm{T}_s:
\Phi\times\Psi\mapsto\mathbb{T}$, $s\geq 0,$ by equality
\eqref{eqn:l2} with $\Phi=\mmf{B}$, $\Psi=\{\phi\}$. From
Lemma~\ref{n:2} and Theorem~\ref{cor:2} it follows that
$$(f'(s)
,\phi)=\lim_{k\to\infty}{(\{\mathrm{F}(t/{g_k})\}^{[g_k\frac{s}{t}]}
\overline{\mmm{Z}}f(0),\phi)},\,\,s\in[0,\, l).$$ Since
$\{n_k\}_{k=1}^{\infty}$ is arbitrary for any fixed
$\phi\in\mmf{X}^{*}$, the Corollary is proved.
\end{proof}

\section{Criteria for closure of operator to be a generator of $C_0$-semigroup.}

\begin{lemma}\label{t:10} Let $\mathrm{F}\in\EuScript{F}_{\mmf{X}}$. Assume that for some $t,\,\, l>0$ and set $\Phi\subset\mmf{X}$
there exists an increasing natural sequence $\{g_k\}_{k=1}^{\infty}$
such that there exists the limit
\begin{equation}
w\mbox{-}\lim\limits_{k
\rightarrow\infty}{\{\mathrm{F}(t/{g_k})\}^{[{g_k}s]} g,}
\label{eqn:m1}
\end{equation}
 for any $g\in \Phi$ and any  rational $s\in[0,\,l)$.
  Then  limit~\eqref{eqn:m1} exists for any $g\in \overline{\span1\{\Phi\}}$, $s\in[0,\,l)$. Furthermore,
  the family of the operators $\mathrm{T}_s: \overline{\span1\{\Phi\}}\mapsto\mmf{X}$,
 $s\in[0,\,l),$ defined by
\begin{equation}
\mathrm{T}_s g=w\mbox{-}\lim\limits_{k
\rightarrow\infty}{\{\mathrm{F}(t/{g_k})\}^{[{g_k}\frac{s}{t}]} g},
\end{equation}
  satisfies the following conditions:

\noindent a) $\mathrm{T}_s g$ is linear w.r.t. $g$  and there exist
$M\geq 1$ and $a\in\mmb{R}$ such that
$$\|\mathrm{T}_s g\|\leq M \exp{(a s)}\|g\|,\,\,s\in[0,\,l).$$

\noindent b) $\mathrm{T}_s g$ is continuous w.r.t. $s\in[0,\,l)$ for
all $g\in \overline{\span1\{\Phi\}}$.

\noindent c) If $f\in
\overline{\span1\{\Phi\}}\cap\mcl{D}(\mmm{F}'(0))$ and
$\mathrm{F}'(0)f\in \overline{\span1\{\Phi\}}$, then there exists
$$(\mathrm{T}_{s} f)^{'}_{s}=\mathrm{T}_{s}\mathrm{F}'(0)f,\,\,s\in[0,\,l).$$

\noindent d) If $\phi\in\mcl{D}({\mmm{F}^*}'(0))$, then
$$(\mathrm{T}_{m}f, \phi)- (\mathrm{T}_{p}f,
\phi)=\int_{p}^{m}{(\mathrm{T}_{s}f, {\mmm{F}^*}'(0)\phi)\,d s},\,\,
m,\,p\in[0,\,l)$$ for any $f\in \overline{\span1\{\Phi\}}$.

\noindent e) If $\phi\in\mcl{D}({\mmm{F}^*}'(0))$, then $
(\mathrm{T}_{s}f, \phi)^{'}_{s}=(\mathrm{T}_{s}f,
{\mmm{F}^*}'(0)\phi),\,s\in[0,\,l),$ for any $f\in
\overline{\span1\{\Phi\}}.$

\noindent f) If $f\in
\overline{\span1\{\Phi\}}\cap\mcl{D}(\mmm{F}'(0))$,
$\mathrm{F}'(0)f\in \overline{\span1\{\Phi\}}$ and
 $\mcl{D}({\mmm{F}^*}'(0))$ is *-dense in $\mmf{X}^*$, then
$\mathrm{T}_{s}f\in\mcl{D}(({\mmm{F}^*}'(0))^*)$, $s\in[0,\,l),$ and
the following statements are valid:

\noindent 1)
$\mcl{D}(({\mmm{F}^*}'(0))^*)\supset\mcl{D}(\mmm{F}'(0))$,
 $({\mmm{F}^*}'(0))^*f=\mathrm{F}'(0)f$.

\noindent 2)
$(\mathrm{T}_{s}f)^{'}_{s}=({\mmm{F}^*}'(0))^*\mathrm{T}_{s}f$ for
all $s\in[0,\,l)$.
\end{lemma}
\begin{proof}[Proof of Lemma~\ref{t:10}.]
The existence of  limit~\eqref{eqn:m1} for any $g\in
\overline{\span1\{\Phi\}}$, $s\geq 0$ can be proved similarly to the
first part of Proposition~\ref{t:1}. Part (a) is trivial. Parts (b)
-- (e) easily follow from Proposition~\ref{t:1}. Let us show part
(f). We can notice that the following chain of the equalities can be
deduced from parts (c) and (e)
 $$((\mmm{T}_s f,
 g))^{'}_s=(\mmm{T}_s\mmm{F}'(0) f, g)=(\mmm{T}_s f,
{\mmm{F}^{*}}'(0) g),\,\,g\in\mcl{D}({\mmm{F}^*}'(0)).$$ Thus, by
*-density of $\mcl{D}({\mmm{F}^*}'(0))$ in $\mmf{X}^*$, we get
$\mmm{T}_s f\in \mcl{D}(({\mmm{F}^{*}}'(0))^{*})$ and
$({\mmm{F}^{*}}'(0))^{*}\mmm{T}_s f=\mmm{T}_s\mmm{F}'(0)f$ for each
$s\geq 0$. From this it follows  that
$f\in\mcl{D}(({\mmm{F}^{*}}'(0))^{*})$ and
$({\mmm{F}^{*}}'(0))^{*}f=\mmm{F}'(0)f$. Hence, by part (c), we get
that $(\mathrm{T}_{s}f)^{'}_{s}=({\mmm{F}^*}'(0))^*\mathrm{T}_{s}f$.
\end{proof}

\begin{theorem}\label{cog:1} Let  $\mathrm{F}\in\EuScript{F}_{\mmf{X}}$. Assume that a linear operator $\mmm{Z}$
has the domain $\mcl{D}(\mmm{Z})\subset\mcl{D}(\mmm{F}'(0))$ and
$$\mmm{Z}f=\mmm{F}'(0)f,\,\,f\in\mcl{D}(\mmm{Z}).$$  Assume also that
$\mmf{A}\subset\mcl{D}(\mmm{Z})$ is a dense linear subset of
$\mmf{X}$ and there exists a fixed $l>0$ such that there exists a
local solution $f:[0,\,l)\mapsto\mcl{D}(\overline{\mmm{Z}})$ of the
system
\begin{eqnarray}
f'(s)&=&\overline{\mmm{Z}}f(s),\,\,s\in[0,\,l),\label{tt:21}\\
f(0)&=&f_0\label{tt:21a}
\end{eqnarray}
for each  $f_0\in\mmf{A}$. Then $\mmm{Z}$ is closable and
$\overline{\mmm{Z}}$ is a generator of $C_0$-semigroup $\mmm{S}$.
Furthermore, the following equality is satisfied:
\begin{eqnarray} \mmm{S}(t) f=\lim_{n\to\infty}{{\mmm{F}(t/n)}^n
f},\, t\geq 0,\label{st:14}
\end{eqnarray} for all $f\in\mmm{X}.$
\end{theorem}
\begin{proof}[Proof of Theorem~\ref{cog:1}.]
By Theorem \ref{cor:2}, there exists the limit
\begin{eqnarray} w\mbox{-}\lim_{n\to\infty}{{\mmm{F}(t/n)}^{[ns/t]}
f},\,\, s\in[0,\,l),\label {st:12}
\end{eqnarray} for each $f\in\mmf{A}$. From the density of
$\mmf{A}$ in $\mmf{X}$ and Lemma \eqref{t:10} we easily infer that
there exists limit~\eqref{st:12} for each $f\in\mmf{X}$. Put
$\mmm{G}(s) f=w\mbox{-}\lim_{n\to\infty}{{\mmm{F}(t/n)}^{[ns/t]}
f},\, s\in[0,\,l),\,f\in\mmf{X}$. From the uniqueness of the local
solution of
 system \eqref{tt:21}--\eqref{tt:21a} we infer that
\begin{eqnarray}
\mmm{G}(s_1)\mmm{G}(s_2)f=\mmm{G}(s_1+s_2)f,\,s_1,\,s_2,\,s_1+s_2\in[0,\,l),\label{st:13}
\end{eqnarray}
for each $f\in\mmf{A}$. From this and part (a) of Lemma \eqref{t:10}
it follows that equality (\ref{st:13}) is valid for each
$f\in\mmf{X}$.
 Define the function
$\mmm{S}:[0,\,\infty)\mapsto\mcl{L}(\mmf{X})$ by the equality
$$\mmm{S}(s)=\{\mmm{G}(l/2)\}^{[2s/l]}\mmm{G}(s-[2s/l]l/2),\,\,s\geq 0.$$ Then, from parts (a), (b) of Lemma \eqref{t:10} and the definition of $\mmm{S}$  we easily
infer that the function $\mmm{S}$ is $C_0$-semigroup. Furthermore,
from part (c) of Lemma \eqref{t:10} it follows that
$\mcl{D}(\mmm{Z})\subset\mcl{D}(\mmm{F}'(0))\subset\mcl{D}(\mmm{S}'(0))$
and, by the closedness of $\mmm{S}'(0)$,
$\overline{\mmm{Z}}f=\mmm{S}'(0)f$ for each
$f\in\mcl{D}(\overline{\mmm{Z}})$. Let us show that
$\mcl{D}(\overline{\mmm{Z}})=\mcl{D}(\mmm{S}'(0))$. Fix
$g\in\mmf{A}$. Then
$$\sum_{k=1}^{n}{{\textstyle\frac{s}{n}\mmm{S}(k\frac{s}{n})g}}\in\mmf{A}$$ and,
by properties of $C_0$-semigroups, there exist the limits
\begin{eqnarray}\lim_{n\to\infty}{\sum_{k=1}^{n}{{\textstyle\frac{s}{n}\mmm{S}(k\frac{s}{n})g}}}&=&\int_{0}^{s}{\mmm{S}(t)g\,d
t},\nonumber\\
\lim_{n\to\infty}{\ov{\mmm{Z}}\sum_{k=1}^{n}{{\textstyle\frac{s}{n}\mmm{S}(k\frac{s}{n})g}}}&=&\int_{0}^{s}{\mmm{S}(t)\ov{\mmm{Z}}g\,d
t}=\mmm{S}(s)g-g,\,\,s\in[0,\,l/2].\nonumber
\end{eqnarray}  Thus, by
the closedness of $\overline{\mmm{Z}}$, we get that
$\int_{0}^{s}{\mmm{S}(t)g\,d t}\in\mcl{D}(\overline{\mmm{Z}})$ and
$$\ov{\mmm{Z}}\int_{0}^{s}{\mmm{S}(t)g\,d t}=\mmm{S}(s)g-g.$$
 Let
$f\in\mmf{X}$. Choose a sequence $\{f_n\}_{n=1}^{\infty},\,
f_n\in\mmf{A},$ such that $\lim\limits_{n\to\infty}{f_n}=f$. Then,
by properties of $C_0$-semigroups, there exist the limits
$$\lim_{n\to\infty}{\int_{0}^{s}{\mmm{S}(t)f_n\,d
t}}=\int_{0}^{s}{\mmm{S}(t)}f\,d t,$$
$$\lim_{n\to\infty}{\overline{\mmm{Z}}\int_{0}^{s}{\mmm{S}(t)f_n\,d
t}}=\lim_{n\to\infty}{\mmm{S}(s)f_n-f_n}=\mmm{S}(s)f-f,\,s\in[0,\,l/2].$$
 Then, by the closedness of $\overline{\mmm{Z}}$,
$$\int_{0}^{s}{\mmm{S}(t)f\,d t}\in\mcl{D}(\overline{\mmm{Z}})$$ and
$$\overline{\mmm{Z}}\int_{0}^{s}{\mmm{S}(t)f\,d t}=\mmm{S}(s)f-f$$
for $s\in[0,\,l/2]$.  Notice that for $f\in\mcl{D}(\mmm{S}'(0))$
there exist the limits
$$\lim_{s\to 0}{s^{-1}\int_{0}^{s}{\mmm{S}(t)f\,d t}}=f,$$
$$\lim_{s\to
0}{s^{-1}\overline{\mmm{Z}}\int_{0}^{s}{\mmm{S}(t)f\,d
t}}=\lim_{s\to 0}{s^{-1}(\mmm{S}(s)f-f)}=\mmm{S}'(0)f.$$ So, by the
closedness of $\overline{\mmm{Z}}$,
$f\in\mcl{D}(\overline{\mmm{Z}})$. Thus, $\overline{\mmm{Z}}$ is a
generator of $C_0$-semigroup $\mmm{S}$ and, by Chernoff's theorem,
we get equality \eqref{st:14}.
\end{proof}

\begin{theorem}\label{main:1} Let  $\mmm{Z}$ be a  densely defined  linear
operator in $\mmf{X}$ and $t>0$. Then  $\mmm{Z}$ is closable  and
its closure is a generator of $C_0$-semigroup if and only if there
exists a function $\mathrm{F}\in\EuScript{F}_{\mmf{X}}$ such that:

\noindent i) $\mcl{D}(\mmm{F}'(0))\supset\mcl{D}(\mmm{Z})$ and
$\mmm{F}'(0)f=\mmm{Z}f,\,\,f\in\mcl{D}(\mmm{Z})$.

\noindent ii) $\mcl{D}({\mmm{F}^*}'(0))$ is *-dense in
$\mathbf{X}^*$.

\noindent iii) There exists a dense linear subspace
$\mmf{A}\subseteq\mcl{D}(\mmm{Z})$ such that for any $f\in\mmf{A}$,
$s\geq 0$ there exists a subsequence $\{f_n^{s}\}_{n=1}^{\infty}$,
$f_n^{s}\in\mcl{D}(\mmm{Z}),$ that satisfies the following
conditions:

\noindent a) $\lim\limits_{n\to\infty}{\|\mathrm{F}^{[n
s]}(\frac{t}{n})f-f_n^{s}\|}=0$.

\noindent b) It is possible to choose a weakly convergent
subsubsequence for any subsequence of the sequence
$\{f_n^{s}\}_{n\in\Nat}$.

\noindent c) It is possible to choose a weakly convergent
subsubsequence for any subsequence of the sequence
$\{\mmm{Z}f_n^{s}\}_{n\in\Nat}$.

\noindent Furthermore, if conditions (i)-(iii) are satisfied, then
$\exp{(s\overline{\mmm{Z}})} f=\lim_{n\to\infty}{{\mmm{F}(s/n)}^n
f},\,s>0$, for all $f\in\mmm{X}.$
\end{theorem}
\begin{proof}[Proof of Theorem~\ref{main:1}.]If $\overline{\mmm{Z}}$ is a generator of $C_0$-semigroup and
$\mathrm{F}(s)=\exp{(s\overline{\mmm{Z}})},\,\,s\geq 0$, then it is
easy to check that conditions (i)-(iii) are satisfied. Indeed, it
follows from  Statement~\ref{st:31} that conditions (i), (ii) are
satisfied. Put $\mmf{A}=\mcl{D}(\mmm{Z})$. We can notice that if
$f\in\mmf{A}$, then
$\exp{(s\overline{\mmm{Z}})}f\in\mcl{D}(\overline{\mmm{Z}})$. From
this we infer condition (iii). Now let us show that conditions
(i)-(iii) imply that $\overline{\mmm{Z}}$ is a generator of
$C_0$-semigroup. Let $\{s_k\}_{k=1}^{\infty}$ be a sequence of all
nonnegative rational numbers. Fix $f\in\mmf{A}$ and an increasing
sequence $\{n_i\}_{i=1}^{\infty}$ of natural numbers. For each
$i\in\mmb{N}$ let the sequence $\{f_n^{s_i}\}_{n=1}^{\infty}$
satisfy condition (iii) for  $s_i$ and $f$. Choose a subsequence
$\{n_i^1\}_{i=1}^{\infty}$ of the sequence $\{n_i\}_{i=1}^{\infty}$
such that there exist the limits
$w\mbox{-}\lim\limits_{i\to\infty}{f_{n_i^1}^{s_1}}$ and
$w\mbox{-}\lim\limits_{i\to\infty}{\mmm{Z}f_{n_i^1}^{s_1}}$. Then,
it follows from part (a) of condition (iii) that there exists the
limit $w\mbox{-}\lim\limits_{i\to\infty}{\mmm{F}^{[{n_i^1}
s_1]}(t/{n_i^1})f}$. Similarly, choose a subsequence $\{n_i^2\}$ of
the sequence $\{n_i^1\}$ such that there exist the limits
  $w\mbox{-}\lim\limits_{i\to\infty}{f_{n_i^2}^{s_2}}$,
  $w\mbox{-}\lim\limits_{i\to\infty}{\mmm{Z}f_{n_i^2}^{s_2}}$,
  $w\mbox{-}\lim\limits_{i\to\infty}{\mmm{F}^{[{n_i^2}
  s_2]}(t/{n_i^2})f}$. In the same way, for any natural $k\geq 3$, we choose a subsequence $\{n_i^k\}_{i=1}^{\infty}$
of the sequence $\{n_i^{k-1}\}_{i=1}^{\infty}$ such that there exist
the limits $w\mbox{-}\lim\limits_{i\to\infty}{f_{n_i^k}^{s_k}}$,
$w\mbox{-}\lim\limits_{i\to\infty}{\mmm{Z}f_{n_i^k}^{s_k}}$,
$w\mbox{-}\lim\limits_{i\to\infty}{\mmm{F}^{[{n_i^k}
s_k]}(t/{n_i^k})f}$. Then we can consider the diagonal sequence
$\{n_i^i\}_{i=1}^{\infty}$ and deduce that there exist the limits
$w\mbox{-}\lim\limits_{i\to\infty}{f_{n_i^i}^{s_k}}$,
$w\mbox{-}\lim\limits_{i\to\infty}{\mmm{Z}f_{n_i^i}^{s_k}}$,
$w\mbox{-}\lim\limits_{i\to\infty}{\mmm{F}^{[n_i^i
s_k]}(t/{n_i^i})f}$ for any $k\in\mmb{N}$. Choose now a sequence
$\{g_l\}_{l=1}^{\infty}$, $g_l\in\mmf{A},$ such that
$\lim\limits_{k\to\infty}{g_k}=\mmm{F}'(0)f$. Choose a subsequence
$d_k$ of the sequence $\{n_k^k\}_{k=1}^{\infty}$ such that
 there exists the limit
  $w\mbox{-}\lim\limits_{l\to\infty}{\mmm{F}^{[d_l
  s_k]}(t/{d_l})g_m}$ for any $k,\,m\in\mmb{N}$. Let
  $$\mmf{B}=\overline{\span1{\{\{f\}\cup\{g_k|k\in\mmb{N}\}\}}}.$$ Then it follows from Lemma~\ref{t:10}
  that there exists the limit $w\mbox{-}\lim\limits_{i\to\infty}{\mmm{F}^{[d_i
  s]}(t/{d_i})g}$ for any $s\geq 0$, $g\in\mmf{B}$. Define
  $$\mmm{T}_s g=w\mbox{-}\lim_{i\to\infty}{\mmm{F}^{[d_i
  s]}(t/{d_i})g}, \,\,g\in\mmf{B}.$$
   From the coincidence of the weak closure of
  $\mmm{Z}$ with the closure in the strong topology it follows that $\mmm{T}_{s_k}f\in\mcl{D}(\overline{\mmm{Z}})$ for any
  $k\in\mmb{N}$ because there exist the limits $w\mbox{-}\lim\limits_{i\to\infty}{f_{d_i}^{s_k}}$,
  $w\mbox{-}\lim\limits_{i\to\infty}{\mmm{Z}f_{d_i}^{s_k}}$.
From part (f)
  of Lemma~\ref{t:10}  and the closedness of the operator $({\mmm{F}^*}'(0))^*$ we infer that
  $$(\mathrm{T}_{s}f)^{'}_{s}=\overline{\mmm{Z}}\mathrm{T}_{s}f$$ for any rational $s\geq 0$.
  Further, we can choose for arbitrary $m>0$ a sequence $\{m_i\}$,
  $m_i\in\mmb{Q}\cap [0,\,\infty)$, such that
  $\lim_{i\to\infty}{m_i}=m$. Then
  $$\lim_{i\to\infty}{\mathrm{T}_{m_i}f}=\mathrm{T}_{m}f$$ and
  $$\lim_{i\to\infty}{\overline{\mmm{Z}}\mathrm{T}_{m_i}f}=
\lim_{i\to\infty}{(\mathrm{T}_{s}f)^{'}_{m_i}f}=\lim_{i\to\infty}{\mathrm{T}_{m_i}\mmm{Z}f}=\mathrm{T}_m\mmm{Z}f.$$
Therefore, by the closability of $\mmm{Z}$,
$\mathrm{T}_{m}f\in\mcl{D}(\overline{\mmm{Z}})$ and
$$(\mathrm{T}_{s}f)^{'}_{s}=\overline{\mmm{Z}}\mathrm{T}_{s}f$$ for
all $s\geq 0$. Consequently, $\mathrm{T}_{s}f$ is a solution of the
equation $$f'(s)=\overline{\mmm{Z}}f(s)$$ with the initial condition
$f(0)=f$. Thus, Theorem~\ref{cog:1} implies that
$\overline{\mmm{Z}}$ is the generator of $C_0$-semigroup $\mmm{S}$
such that
$$\mmm{S}(s) f=\lim_{n\to\infty}{{\mmm{F}(s/n)}^n
f},\, s\geq 0,$$ for each $f\in\mmm{X}.$
\end{proof}

\begin{theorem}\label{t:4} Let  $\mmm{Z}$ be a densely defined linear
operator in $\mmf{X}$ and $t>0$. Assume also that there exists a
function $\mathrm{F}\in\EuScript{F}_{\mmf{X}}$ such that:

\noindent i) $\mcl{D}(\mmm{F}'(0))\supset\mcl{D}(\mmm{Z})$  and
$\mmm{F}'(0)f=\mmm{Z}f,\,\,f\in\mcl{D}(\mmm{Z})$.

\noindent ii) $\mcl{D}({\mmm{F}^*}'(0))$  *-dense in $\mathbf{X}^*$.

\noindent Then operator $\mmm{Z}$ is closable  and its closure is a
generator of $C_0$-semigroup iff there exists a dense linear
subspace $\mmf{A}\subseteq\mcl{D}(\mmm{Z})$ such that for all
$f\in\mmf{A}$, $s\geq 0$ there exists a sequence
$\{f_n^{s}\}_{n=1}^{\infty}$, $f_n^{s}\in\mcl{D}(\mmm{Z}),$
satisfying the  following conditions:

\noindent a) $\lim\limits_{n\to\infty}{\|\mathrm{F}^{[n
s]}(\frac{t}{n})f-f_n^{s}\|}=0$.

\noindent b) It is possible to choose a weakly convergent
subsubsequence for any subsequence of the sequence
$\{\mmm{Z}f_n^{s}\}_{n\in\Nat}$.

\noindent c) It is possible to choose a weakly convergent
subsubsequence for any subsequence of the sequence
$\{f_n^{s}\}_{n\in\Nat}$.
\end{theorem}
\begin{proof}[Proof of Theorem~\ref{t:4}.]
\begin{trivlist}
\item{($\Rightarrow$)} If operator $\mmm{Z}$ is closable  and its closure is a
generator of $C_0$-semigroup, then, by Chernoff's Theorem, it is
enough to put $\mmf{A}=\mcl{D}(\mmm{Z})$ and  choose sequences
$\{f_n^{s}\}_{n=0}^{\infty},\,\,f_n^{s}\in\mmf{A},\,\,s\geq 0,$ such
that there exist the limits
$\lim\limits_{n\to\infty}{f_n^{s}}=\exp{(s\ov{\mmm{Z}})}f$,
$\lim\limits_{n\to\infty}{\mmm{Z}f_n^{s}}=\ov{\mmm{Z}}\exp{(s\ov{\mmm{Z}})}f$.
\item{($\Leftarrow$)} It follows from Theorem~\ref{main:1}.
\end{trivlist}
\end{proof}

\begin{corollary}\label{t:3} Let $t>0$. Assume that $\mmm{C}$ and
$\mmm{D}$ are generators of $C_0$-semigroups $\exp{(s\mmm{C})}$ and
$\exp{(s\mmm{D})}$ in $\mmf{X}$. Assume also that there exists a set
$\mmf{B}\subset\mcl{D}((\exp(s\mmm{C})^*)'_{s=0})\cap\mcl{D}((\exp(s\mmm{D})^*)'_{s=0})$
and the following conditions are satisfied:

\noindent i) $\mcl{D}(\mmm{C})\cap\mcl{D}(\mmm{D})$ is dense in
$\mmf{X}$.

\noindent ii) $\mmf{B}$ is *-dense in $\mmf{X^*}$.

\noindent iii) There exist $a\in\mathbb{R}$ and $M\geq 1$ such that
$$\|\{\exp{{\textstyle(\frac{s}{n}\mmm{C})}}\exp{{\textstyle(\frac{s}{n}\mmm{D})}}\}^m\|\leq
M\exp{({\textstyle a s\frac{m}{n}})}$$ for all $n,\, m\in\mmb{N}$,
$s>0$.

\noindent iv) Function $g(s,x)=\exp{(s\mmm{D}^*)}x,\,\, s\geq 0,\,\,
x\in \mmf{X},$ is continuous at  $s=0$ for $x\in
(\exp(s\mmm{C})^*)'_{s=0}(\mmf{B})$.

\noindent Then the sum of $\mmm{C}$ and $\mmm{D}$ is a closable
operator
 and its closure is a generator of $C_0$-semigroup iff there exists
 a
dense linear subspace
$\mmf{A}\subseteq\mcl{D}(\mmm{C})\cap\mcl{D}(\mmm{D})$ such that for
all  $f\in\mmf{A}$, $s\geq 0$ there exists a sequence
$\{f_n^{s}\}_{n=1}^{\infty}$,
$f_n^{s}\in\mcl{D}(\mmm{C})\cap\mcl{D}(\mmm{D}),$ satisfying the
following conditions:

\noindent a) $\lim\limits_{n\to\infty}{\|(\exp{{(\textstyle
\frac{t}{n}\mmm{C})}}\exp{{(\textstyle\frac{t}{n}\mmm{D})}})^{[n
s]}f-f_n^{s}\|}=0$.

\noindent b) It is possible to choose a weakly convergent
subsubsequence for any subsequence of the sequence
$\{f_n^{s}\}_{n=1}^{\infty}$.

\noindent c) It is possible to choose a weakly convergent
subsubsequence for any subsequence of the sequence
$\{(\mmm{C}+\mmm{D})f_n^{s}\}_{n=1}^{\infty}$.

\noindent Furthermore, if conditions (a)-(c) are satisfied, then
$$\exp{(s(\overline{\mmm{C}+\mmm{D}}))} f=\lim_{n\to\infty}{\{\exp{{\textstyle (\frac{s}{n}\mmm{C})}}\exp{{\textstyle (\frac{s}{n}\mmm{D})}}\}^n
f},\, s\geq 0,$$ for all $f\in\mmm{X}.$
\end{corollary}
\begin{proof}[Proof of Corollary~\ref{t:3}.] Let $\mathrm{F}(s)=\exp(s \mmm{C})\exp(s \mmm{D})$.
Since
\begin{eqnarray}
 & &\lim\limits_{s \rightarrow 0}{s^{-1}(\mmm{F}(s)f-\mmm{F}(0)f)}=\lim\limits_{s \rightarrow
0}{s^{-1}(\exp(s \mmm{C})\exp(s \mmm{D})f-\exp(s
\mmm{C})f)}\nonumber\\&+& \lim\limits_{s \rightarrow
0}{s^{-1}(\exp(s
\mmm{C})f-f)}=\mmm{C}f+\mmm{D}f,\,\,f\in\mcl{D}(\mmm{C})\cap\mcl{D}(\mmm{D}),\nonumber
\end{eqnarray}
and
\begin{eqnarray} \lim\limits_{s \rightarrow
0}{s^{-1}(\mmm{F}^*(s)-\mmm{F}(0))f}=\lim\limits_{s \rightarrow
0}{s^{-1}\{(\exp(s \mmm{D}))^*(\exp(s \mmm{C}))^*-(\exp(s
\mmm{D}))^*\}f}\nonumber\\+ \lim\limits_{s \rightarrow
0}{s^{-1}\{(\exp(s \mmm{D}))^*f-f\}}=\lim\limits_{l \rightarrow
0}{(\exp(l
\mmm{D}))^*(\exp(s\mmm{C})^*)'|_{s=0}f}\nonumber\\+(\exp(s\mmm{D})^*)'|_{s=0}f=
(\exp(s\mmm{C})^*)'|_{s=0}f+(\exp(s\mmm{D})^*)'|_{s=0}f,\,\,f\in
\mmf{B},\nonumber
\end{eqnarray}
we infer that $\mcl{D}(\mmm{F}'(0))$ is dense in $\mmm{X}$ and
$\mcl{D}({\mmm{F}^{*}}'(0))$ is *-dense in $\mmm{X}^*$. Therefore,
we can apply Theorem~\ref{t:4} with $\mmm{Z}=\mmm{C}+\mmm{D}$.
\end{proof}

\begin{lemma}\label{t:44} Let $\mathrm{F}\in\EuScript{F}_{\mmf{X}}$. Assume that for some $t>0$ and set $\Phi\subset\mmf{X}$
there exists an increasing natural sequence $\{g_k\}_{k=1}^{\infty}$
such that there exists the limit
\begin{equation}
\lim\limits_{k \rightarrow\infty}{\{\mathrm{F}(t/{g_k})\}^{[{g_k}s]}
g} \label{eqn:v1}
\end{equation}
 for any $g\in \Phi$ and any nonnegative rational $s$.
  Then  limit~\eqref{eqn:v1} exists for any $g\in \overline{\span1\{\Phi\}}$, $s\geq 0$. Furthermore,
  the family of the operators $\mathrm{T}_s: \overline{\span1\{\Phi\}}\mapsto\mmf{X}$,
 $s\geq 0,$ defined by
\begin{equation}
\mathrm{T}_s g=\lim\limits_{k
\rightarrow\infty}{\{\mathrm{F}(t/{g_k})\}^{[{g_k}\frac{s}{t}]}
g},\nonumber
\end{equation}
  satisfies the following conditions:

\noindent a) If $\mathrm{T}_{s} f\in \overline{\span1\{\Phi\}}$ for
some $s\geq 0$, then
\begin{eqnarray} \mathrm{T}_{s+l}f
=\mathrm{T}_{l}\mathrm{T}_{s}f,\,\,l> 0.\nonumber
\end{eqnarray}

\noindent b) Let $f\in
\overline{\span1\{\Phi\}}\cap\mcl{D}(\ov{\mmm{F}'(0)})$ and
$\ov{\mathrm{F}'(0)}f\in \overline{\span1\{\Phi\}}$. Assume also
that there exists a sequence
$\{f_n\}_{n=1}^{\infty},\,f_n\in\overline{\span1\{\Phi\}}\cap\mcl{D}(\mmm{F}'(0))$,
such that $\lim_{n\to\infty}{f_n}=f$ and
$\lim_{n\to\infty}{\mmm{F}'(0)f_n}=\ov{\mmm{F}'(0)}f$. Then there
exists
$$(\mathrm{T}_{s} f)^{'}_{s}=\mathrm{T}_{s}\ov{\mathrm{F}'(0)}f,\,\,s\geq
0.$$

\noindent c) Let $\mathrm{T}_{s} f\in
\overline{\span1\{\Phi\}}\cap\mcl{D}(\ov{\mmm{F}'(0)})$, $ f\in
\overline{\span1\{\Phi\}}\cap\mcl{D}(\mmm{F}'(0))$
 and $\ov{\mmm{F}'(0)}\mathrm{T}_{s}
f,\,\mmm{F}'(0)f\in\overline{\span1\{\Phi\}}$ for some $s\geq 0$.
Assume also that there exists a sequence
$\{f_n\}_{n=1}^{\infty},\,f_n\in\overline{\span1\{\Phi\}}\cap\mcl{D}(\mmm{F}'(0))$,
such that $\lim_{n\to\infty}{f_n}=\mathrm{T}_{s} f$ and
$\lim_{n\to\infty}{\mmm{F}'(0)f_n}=\ov{\mmm{F}'(0)}\mathrm{T}_{s}
f$. Then there exists
\begin{eqnarray}
(\mathrm{T}_{s} f)^{'}_{s}=
\overline{\mathrm{F}'(0)}\mathrm{T}_{s}f,\,\,s\geq 0.\nonumber
\end{eqnarray}
\end{lemma}
\begin{proof}[Proof of Lemma~\ref{t:44}.]
The existence of  limit~\eqref{eqn:v1} for any $g\in
\overline{\span1\{\Phi\}}$, $s\geq 0$ can be proved similarly to the
first part of Proposition~\ref{t:1}.
\begin{trivlist}
\item{(a)} Part (a) follows from the chain of the equalities:
\begin{eqnarray}
 \mathrm{T}_l \mathrm{T}_s f=\lim\limits_{k
\rightarrow\infty}{\{\mathrm{F}(t/{g_k})\}^{[{g_k}l]} \mathrm{T}_s
f}=\lim\limits_{k
\rightarrow\infty}{\{\mathrm{F}(t/{g_k})\}^{[{g_k}l]}
(\mathrm{T}_s-\{\mathrm{F}(t/{g_k})\}^{[{g_k}s]}) f}\nonumber\\+
\lim\limits_{k \rightarrow\infty}{\{\mathrm{F}(t/{g_k})\}^{[{g_k}l]}
\{\mathrm{F}(t/{g_k})\}^{[{g_k}s]} f}=\lim\limits_{k
\rightarrow\infty}{\{\mathrm{F}(t/{g_k})\}^{[{g_k}l]+[{g_k}s]}
f}=\mathrm{T}_{l+s} f\nonumber
\end{eqnarray}
\item{(b)} By part (c) of Lemma \ref{t:10}, we have
\begin{eqnarray}
\mmm{T}_s f-\mmm{T}_{r}f=\lim_{n\to\infty}{\mmm{T}_s
f_n-\mmm{T}_{r}f_n}=\lim_{n\to\infty}{\int_{r}^{s}{\mmm{T}_l
\mathrm{F}'(0)f_n \,d l}}=\int_{r}^{s}{\mmm{T}_l
\ov{\mathrm{F}'(0)}f \,d l}\nonumber
\end{eqnarray}
for any  $r,\,\,s\geq 0$. From this we easily infer part (b).

\item{(c)} From part (c) of Lemma \ref{t:10} and parts (a), (b) of Lemma \ref{t:44} it follows that
\begin{eqnarray}
(\mathrm{T}_{s} f)^{'}_{s}= \mathrm{T}_{s}\mathrm{F}'(0)f=\lim_{l\to
0}{l^{-1}(\mathrm{T}_{s+l}f-\mathrm{T}_{s}f)}=\lim_{l\to
0}{l^{-1}(\mathrm{T}_{l}-\mathrm{I})\mathrm{T}_s
f}=\overline{\mathrm{F}'(0)}\mathrm{T}_{s}f.\nonumber
\end{eqnarray}
\end{trivlist}

\end{proof}

The  following two results are respectively extensions of
Corollaries 4, 5 from \cite{Nek}.

\begin{theorem}\label{t:54} Let  $\mmm{Z}$ be a densely defined linear operator in $\mmf{X}$
and $t>0$. Assume also that there exists a function
$\mathrm{F}\in\EuScript{F}_{\mmf{X}}$ such that
$\mcl{D}(\mmm{F}'(0))\supset\mcl{D}(\mmm{Z})$  and
$\mmm{F}'(0)f=\mmm{Z}f,\,\,f\in\mcl{D}(\mmm{Z})$. Then operator
$\mmm{Z}$ is closable  and its closure is a generator of
$C_0$-semigroup iff there exists a dense linear subspace
$\mmf{A}\subseteq\mcl{D}(\mmm{Z})$ such that for all $f\in\mmf{A}$,
$s\geq 0$ there exists a sequence $\{f_n^{s}\}_{n=1}^{\infty}$,
$f_n^{s}\in\mcl{D}(\mmm{Z}),$ satisfying the  following conditions:

\noindent a) $\lim\limits_{n\to\infty}{\|\mathrm{F}^{[n
s]}(\frac{t}{n})f-f_n^{s}\|}=0$.

\noindent b) The sets $\{\mmm{Z}f_n^{s}\}_{n=1}^{\infty}$ and
$\{f_n^{s}\}_{n=1}^{\infty}$ are precompact.

\end{theorem}

\begin{proof}[Proof of Theorem~\ref{t:54}.]
\begin{trivlist}
\item{($\Rightarrow$)} If operator $\mmm{Z}$ is closable  and its closure is a
generator of $C_0$-semigroup, then, by Chernoff's Theorem, it is
enough to put $\mmf{A}=\mcl{D}(\mmm{Z})$ and  choose sequences
$\{f_n^{s}\}_{n=0}^{\infty},\,\,f_n^{s}\in\mmf{A},\,\,s\geq 0,$ such
that there exist the limits
$\lim\limits_{n\to\infty}{f_n^{s}}=\exp{(s\ov{\mmm{Z}})}f$,
$\lim\limits_{n\to\infty}{\mmm{Z}f_n^{s}}=\ov{\mmm{Z}}\exp{(s\ov{\mmm{Z}})}f$.
\item{($\Leftarrow$)}
From condition (a) it follows that $\{\mathrm{F}^{[n
s]}(\frac{t}{n})f\}_{n=1}^{\infty}$ is precompact for each
$f\in\mmf{A}$ and $s>0$. Then, by the density $\mmf{A}$ in
$\mmf{X}$, $\{\mathrm{F}^{[n s]}(\frac{t}{n})f\}_{n=1}^{\infty}$ is
precompact for each $f\in\mmf{X}$ and $s>0$. Fix $f\in\mmf{A}$. As
in the proof of Theorem \ref{main:1}, we can choose an increasing
sequence $\{g_n\}_{n=1}^{\infty}$ such that  there exist the limits
$\lim\limits_{n\to\infty}{\mathrm{F}^{[g_n s]}(\frac{t}{g_n})f}$,
 $\lim\limits_{n\to\infty}{\mathrm{Z}f_{g_n}^{s}}$,
$\lim\limits_{n\to\infty}{\mathrm{F}^{[g_n
l]}(\frac{t}{g_n})f_{k}^{s}}$ and
$\lim\limits_{n\to\infty}{\mathrm{F}^{[g_n
l]}(\frac{t}{g_n})\mathrm{Z}f_{k}^{s}}$ for each
$s,\,l\in[0,\,\infty)\cap\mmb{Q}$ and $k\in\mmb{N}$. Let
  $$\mmf{B}=\overline{\span1{\{\{f\}\cup\{\mmm{Z}f_{k}^{s}|k\in\mmb{N},\,s\in[0,\,\infty)\cap\mmb{Q}\}\cup\{f_{k}^{s}|k\in\mmb{N},\,s\in[0,\,\infty)\cap\mmb{Q}\}\}}}.$$
   Then it follows from Lemma~\ref{t:44}
  that there exists the limit $\lim\limits_{i\to\infty}{\mmm{F}^{[d_i
  s]}(t/{d_i})g}$ for any $s\geq 0$, $g\in\mmf{B}$. Define
  $$\mmm{T}_s g=\lim_{i\to\infty}{\mmm{F}^{[d_i
  s]}(t/{d_i})g}, \,\,g\in\mmf{B}.$$ By condition (a), there exist the limits $\lim\limits_{i\to\infty}{f_{g_n}^{s}}=\mmm{T}_s f$ and
  $\lim\limits_{i\to\infty}{\mmm{Z}f_{g_n}^{s}}$ for each $s\in[0,\,\infty)\cap\mmb{Q}$. So, by the closability of
  $\mmm{Z}$,
  $\mmm{T}_{s}f\in\mcl{D}(\overline{\mmm{Z}})\cap\mmf{B}$
   and $$\lim_{i\to\infty}{\mmm{Z}f_{g_n}^{s}}= \ov{\mmm{F}'(0)}\mmm{T}_{s}f\in\mmf{B}$$ for any
  $s\in[0,\,\infty)\cap\mmb{Q}$. Thus, by part (c) of Lemma~\ref{t:44}, we get
  $$(\mathrm{T}_{s}f)^{'}_{s}=\overline{\mmm{Z}}\mathrm{T}_{s}f$$ for any rational $s\geq 0$.
  Choose for arbitrary $m>0$ a sequence $\{m_i\}$,
  $m_i\in\mmb{Q}\cap [0,\,\infty)$, such that
  $\lim\limits_{i\to\infty}{m_i}=m$. Then
  $$\lim_{i\to\infty}{\mathrm{T}_{m_i}f}=\mathrm{T}_{m}f$$ and
  $$\lim_{i\to\infty}{\overline{\mmm{Z}}\mathrm{T}_{m_i}f}=
\lim_{i\to\infty}{(\mathrm{T}_{s}f)^{'}_{m_i}f}=\lim_{i\to\infty}{\mathrm{T}_{m_i}\mmm{Z}f}=\mathrm{T}_m\mmm{Z}f.$$
Therefore, by the closability of $\mmm{Z}$,
$\mathrm{T}_{m}f\in\mcl{D}(\overline{\mmm{Z}})$ and
$$(\mathrm{T}_{s}f)^{'}_{s}=\overline{\mmm{Z}}\mathrm{T}_{s}f$$ for
all $s\geq 0$. Consequently, $\mathrm{T}_{s}f$ is a solution of the
equation $$f'(s)=\overline{\mmm{Z}}f(s)$$ with the initial condition
$f(0)=f$. Thus, Theorem~\ref{cog:1} implies that
$\overline{\mmm{Z}}$ is a generator of $C_0$-semigroup
 and, by Chernoff's theorem,
$$\exp{(s\ov{\mmm{Z}})} f=\lim_{n\to\infty}{{\mmm{F}(s/n)}^n
f},\,\,f\in\mmm{X}.$$
\end{trivlist}
\end{proof}
From Theorem~\ref{t:54} we get the following result.
\begin{corollary} \label{t:55} Let $t>0$. Assume that $\mmm{C}$ and
$\mmm{D}$ are generators of $C_0$-semigroups $\exp{(s\mmm{C})}$ and
$\exp{(s\mmm{D})}$. Assume also that the following conditions are
satisfied:

\noindent i) $\mcl{D}(\mmm{C})\cap\mcl{D}(\mmm{D})$ is dense in
$\mmf{X}$.

\noindent ii) There exists $a\in\mathbb{R}$ and $M\geq 1$ such that
$$\|\{\exp{{\textstyle(\frac{s}{n}\mmm{C})}}\exp{{\textstyle(\frac{s}{n}\mmm{D})}}\}^m\|\leq
M\exp{({\textstyle a s\frac{m}{n}})}$$ for all $n,\, m\in\mmb{N}$,
$s\in\mmb{R}$.

\noindent Then the sum of $\mmm{C}$ and $\mmm{D}$ is a closable
operator
 and its closure is a generator of $C_0$-semigroup iff there exists
 a
dense linear subspace
$\mmf{A}\subseteq\mcl{D}(\mmm{C})\cap\mcl{D}(\mmm{D})$ such that for
all  $f\in\mmf{A}$, $s\geq 0$ there exists a sequence
$\{f_n^{s}\}_{n=1}^{\infty}$,
$f_n^{s}\in\mcl{D}(\mmm{C})\cap\mcl{D}(\mmm{D}),$ satisfying the
following conditions:

\noindent a) $\lim\limits_{n\to\infty}{\|(\exp{{(\textstyle
\frac{t}{n}\mmm{C})}}\exp{{(\textstyle\frac{t}{n}\mmm{D})}})^{[n
s]}f-f_n^{s}\|}=0$.

\noindent b) The sets $\{f_n^{s}\}_{n=1}^{\infty}$ and
$\{(\mmm{C}+\mmm{D})f_n^{s}\}_{n=1}^{\infty}$ are precompact.

\noindent Furthermore, if conditions (a)-(b) are satisfied, then
$$\exp{(s(\overline{\mmm{C}+\mmm{D}}))} f=\lim_{n\to\infty}{\{\exp{{\textstyle (\frac{s}{n}\mmm{C})}}\exp{{\textstyle (\frac{s}{n}\mmm{D})}}\}^n
f},\, s>0,$$ for all $f\in\mmm{X}.$
\end{corollary}

If $\mmf{X}$ is a reflexive separable Banach space or a Hilbert
space,
 then the formulations of
 Theorems~\ref{main:1},~\ref{t:4} and Corollary~\ref{t:3} can be
considerably simplified and, by Statement \ref{st:32}, we get the
following Corollaries of the results mentioned above:
\begin{corollary}\label{news:3}Let  $\mmm{Z}:\mcl{B}\supset \mcl{D}(\mmm{Z})\to\mcl{B}$ be a linear operator, $\mcl{D}(\mmm{Z})$ be dense in $\mcl{B}$
and $t>0$. Operator $\mmm{Z}$ is closable  and its closure is a
generator
 of $C_0$-semigroup iff there exists
$\mathrm{F}:[0,\infty)\mapsto\mathcal{L}(\mcl{B})$ such that:

\noindent i) $\mathrm{F}(0)=\mathrm{I}$ and there exists
$a\in\mathbb{R}$, $M\geq 1$ such that
$\|\mathrm{F}^{m}(\frac{s}{n})\|\leq M\exp(a s\frac{m}{n})$ for all
$n,\, m\in\mmb{N}$ and $s\geq 0$.

\noindent ii) $\mcl{D}(\mmm{F}'(0))\supset\mcl{D}(\mmm{Z})$ and
$\mmm{F}'(0)f=\mmm{Z}f,\,при\,f\in\mcl{D}(\mmm{Z})$.

\noindent iii) $\mcl{D}({\mmm{F}^*}'(0))$ is *-dense in $\mcl{B}^*$.

\noindent iv) There exists a dense linear subspace
$\mmf{A}\subseteq\mcl{D}(\mmm{Z})$ such that for any $f\in\mmf{A}$,
$s\geq 0$ there exists a sequence $\{f_n\}_{n=1}^{\infty}$,
$f_n\in\mcl{D}(\mmm{Z})$ satisfying the following identity
$$\lim_{n\to\infty}{\|\mathrm{F}^{[n
s]}({\textstyle\frac{t}{n}})f-f_n\|}=0$$ and, furthermore,
$\{\mmm{Z}f_n\}_{n=1}^{\infty}$ is a bounded sequence.

\noindent Furthermore, if conditions (i)-(iv) are satisfied, then
$\exp{(s\overline{\mmm{Z}})}
f=\lim\limits_{n\to\infty}{{\mmm{F}(s/n)}^n f}, \,s>0,$ for all
$f\in\mcl{B}.$
\end{corollary}

\begin{corollary}\label{new:1}Let  $\mmm{Z}:\mcl{B}\supset \mcl{D}(\mmm{Z})\to\mcl{B}$ be a linear operator, $\mcl{D}(\mmm{Z})$ be dense in $\mcl{B}$
and $t>0$. Assume that there exists
$\mathrm{F}:[0,\infty)\mapsto\mathcal{L}(\mcl{B})$ such that the
following conditions are satisfied:

\noindent i) $\mathrm{F}(0)=\mathrm{I}$ and there exists
$a\in\mathbb{R}$, $M\geq 1$ such that
$$\|\mathrm{F}^{m}({\textstyle\frac{s}{n}})\|\leq M\exp{({\textstyle a s\frac{m}{n}})}$$ for
all $n,\, m\in\mmb{N}$ and $s\geq 0$.

\noindent ii) $\mcl{D}(\mmm{F}'(0))\supset\mcl{D}(\mmm{Z})$ and
$\mmm{F}'(0)f=\mmm{Z}f,\,\,f\in\mcl{D}(\mmm{Z})$.

\noindent iii) $\mcl{D}({\mmm{F}^*}'(0))$ is *-dense in
$\mathbf{B}^*$.

\noindent Then operator $\mmm{Z}$ is a closable operator and its
closure is a generator
 of $C_0$-semigroup iff there exists
a dense linear subspace $\mmf{A}\subseteq\mcl{D}(\mmm{Z})$ such that
for any $f\in\mmf{A}$, $s\geq 0$ there exists a sequence
$\{f_n\}_{n=1}^{\infty}$, $f_n\in\mcl{D}(\mmm{Z}),$ satisfying the
following identity $\lim\limits_{n\to\infty}{\|\mathrm{F}^{[n
s]}({\textstyle\frac{t}{n}})f-f_n\|}=0$ and, furthermore,
$\{\mmm{F}'(0)f_n\}_{n=1}^{\infty}$ is a bounded sequence.
\end{corollary}

\begin{corollary}\label{corol} Let $t>0$. Suppose that $\mmm{C}$ and
$\mmm{D}$ are generators of $C_0$-semigroups $\exp(s\mmm{C})$ and
$\exp(s\mmm{D})$ in $\mcl{B}$. Assume also that the following
conditions are satisfied:

\noindent i) $\mcl{D}(\mmm{C})\cap\mcl{D}(\mmm{D})$ is dense in
$\mcl{B}$.

\noindent ii) $\mcl{D}(\mmm{C}^*)\cap\mcl{D}(\mmm{D}^*)$ is dense in
$\mcl{B^*}$.

\noindent iii) There exists $a\in\mathbb{R}$, $M\geq 1$ such that
the following inequality holds
$$\|\{\exp{{\textstyle (\frac{s}{n}\mmm{C})}}\exp{{\textstyle (\frac{s}{n}\mmm{D})}}\}^m\|\leq M\exp{({\textstyle a
s\frac{m}{n}})}$$ for all $n,\, m\in\mmb{N}$, $s>0$.

\noindent Then the sum of operators $\mmm{C}$ and $\mmm{D}$ is a
closable operator and its closure is a generator of $C_0$-semigroup
iff there exists a dense linear subspace
$\mmf{A}\subseteq\mcl{D}(\mmm{C})\cap\mcl{D}(\mmm{D})$ such that for
any $f\in\mmf{A}$, $s\geq 0$ there exists a sequence
$\{f_n\}_{n=1}^{\infty}$,
$f_n\in\mcl{D}(\mmm{C})\cap\mcl{D}(\mmm{D}),$ satisfying the
following conditions:

\noindent a) $\lim\limits_{n\to\infty}{\|(\exp{({\textstyle
\frac{t}{n}\mmm{C}})}\exp{({\textstyle\frac{t}{n}\mmm{D}})})^{[n
s]}f-f_n\|}=0.$

\noindent b) $\{(\mmm{C}+\mmm{D})f_n\}_{n=0}^{\infty}$ is a bounded
sequence.

\noindent Furthermore, if conditions (a)-(b) are satisfied, then
$$\exp{(s(\overline{\mmm{C}+\mmm{D}}))} f=\lim_{n\to\infty}{\{\exp{{\textstyle (\frac{s}{n}\mmm{C})}}\exp{{\textstyle (\frac{s}{n}\mmm{D})}}\}^n
f},\, \, s>0,$$ for all $f\in\mcl{B}.$
\end{corollary}
\begin{remark} If operators $i\mmm{C}$ and $i\mmm{D}$ are self-adjoint, then conditions (ii), (iii) in Corollary \ref{corol} can be omitted.
\end{remark}

\end{document}